\documentclass[12pt]{amsart}
\usepackage{amsmath,amsthm,amsfonts,amssymb}
\begin{document}

\newcommand{\A}{{\mathbb A}}
\newcommand{\B}{{\mathbb B}}
\newcommand{\C}{{\mathbb C}}
\renewcommand{\H}{{\mathbb H}}
\newcommand{\N}{{\mathbb N}}
\newcommand{\M}{{\mathcal M}}
\newcommand{\Q}{{\mathbb Q}}
\newcommand{\Ok}{{\mathcal O}_k}
\newcommand{\Z}{{\mathbb Z}}
\renewcommand{\O}{{\mathcal O}}
\renewcommand{\P}{{\mathbb P}}
\newcommand{\R}{{\mathbb R}}
\newcommand{\rc}{\subset}
\newcommand{\rank}{\mathop{rank}}
\newcommand{\trace}{\mathop{tr}}
\newcommand{\dimc}{\mathop{dim}_{\C}}
\newcommand{\Lie}{\mathop{Lie}}
\newcommand{\Hom}{\mathop{Hom}}
\newcommand{\Spec}{\mathop{Spec}}
\newcommand{\Aut}{\mathop{{\rm Aut}}}
\newcommand{\Auto}{\mathop{{\rm Aut}_{\mathcal O}}}
\newcommand{\alg}[1]{{\mathbf #1}}

\newtheorem{lemma}{Lemma}[section]

\newtheorem{definition}[lemma]{Definition}
\let\dfn\definition

\newtheorem{corollary}[lemma]{Corollary}

\newtheorem{example}{Example}
\newtheorem*{remark}{Remark}
\newtheorem*{observation}{Observation}
\newtheorem{fact}[lemma]{Fact}
\newtheorem*{remarks}{Remarks}

\newtheorem{proposition}[lemma]{Proposition}
\newtheorem{theorem}[lemma]{Theorem}

\numberwithin{equation}{section}
\def\labelenumi{\rm(\roman{enumi})}

\title{
Tame Discrete Sets in Algebraic Groups
}

\author {J\"org Winkelmann}

\begin{abstract}
\end{abstract}
\subjclass{}

\subjclass{32M17}

\keywords{tame discrete sets, linear algebraic group, Stein manifold}

\address{
J\"org Winkelmann \\
Lehrstuhl Analysis II \\
Mathematisches Institut \\
Ruhr-Universit\"at Bochum\\
44780 Bochum \\
Germany\\
}
\email{joerg.winkelmann@rub.de
\\  
ORCID: 0000-0002-1781-5842
}

\maketitle

\section{Introduction}

For discrete subsets in $\C^n$ the notion of being ``tame'' was defined
in the important paper of Rosay and Rudin \cite{RR}. 
A discrete subset $D\subset\C^n$ is called ``tame''
if and only if there exists an automorphism $\phi$ of $\C^n$ such that
$\phi(D)=\N\times \{0\}^{n-1}$.
(In this paper a subset $D$ of a topological space $X$ is called 
a ``discrete subset'' if every point $p$ in $X$ admits an open neighbourhood
$W$ such that $W\cap D$ is finite.)

In \cite{JW-TAME1} we introduced a new definition of tameness which
applies to arbitrary complex manifolds and agrees with the tameness
notion of Rosay and Rudin for the case $X\simeq\C^n$.

\begin{dfn}\label{def-tame}
Let $X$ be a complex manifold. 
An infinite discrete subset $D$ is called (weakly) {\em tame} if for every
exhaustion function $\rho:X\to\R^+$ and every map $\zeta:D\to\R^+$
there exists an automorphism $\phi$ of $X$ such that
$\rho(\phi(x))\ge \zeta(x)$ for all $x\in D$.
\end{dfn}

In terms of Nevanlinna theory a similar tameness notion
may be formulated as follows:
A discrete set is tame if its counting function $N(r,D)$ may be made
as small as desired. (See \S\ref{sect-nevanlinna} for details.)

First results about this notion obtained in \cite{JW-TAME1}
suggested that best results are to be expected in the case
of complex manifolds whose automorphism group is very large in a certain
sense. In particular, in \cite{JW-TAME1} we proved some significant
results for the case $X=SL_n(\C)$ and found some evidence suggesting
that our tameness
notion might not work very well for non-Stein manifolds
like $\C^n\setminus\{(0,\ldots,0)\}$ or partially hyperbolic
complex manifolds like $\Delta\times\C$.

In this paper we show that for complex manifolds biholomorphic to
complex linear algebraic groups without non-trivial morphism to the
multiplicative group $\C^*$ we obtain a theory of tame discrete
sets essentially
as strong as the theory which Rosay and Rudin developed
for $\C^n$.

Andrist and Ugolini 
(\cite{AU}) have proposed a different notion, namely the
following:

\begin{dfn}\label{def-strong-tame}
Let $X$ be a complex manifold. 
An infinite discrete subset $D$ is called (strongly) {\em tame} 
if for every injective map $f:D\to D$ there exists
an automorphism $\phi$ of $X$ such that $\phi(x)=f(x)$ 
for all $x\in D$.
\end{dfn}

They proved that every complex linear algebraic groups admits
a {\em strongly tame} discrete subset (\cite{AU}).

It is easily verified that ``strongly tame'' implies ``weakly tame''
(see \cite{JW-TAME1}). For $X\simeq\C^n$
and $X\simeq SL_n(\C)$ both tameness notions coincide
(\cite{JW-TAME1}). Furthermore,
for $X=\C^n$ both notions agree with tameness as defined by Rosay and 
Rudin.

However, for arbitrary manifolds ``strongly tame'' and ``weakly tame''
are not equivalent (\cite{JW-TAME1}).

In this article, unless explicitly stated otherwise, tame always
means weakly tame, i.e., tame in the sense of Definition~\ref{def-tame}.

\section{Main results}
\subsection{Equivalence of tameness notions}

This article is focused on the study of tame discrete subset
in complex linear-algebraic
groups whose character group is trivial, i.e., which do not admit
any non-constant morphism of algebraic groups to the multiplicative
group $\C^*$. For simplicity, we call these algebraic groups
``character-free''.

First we show that for this class of complex manifolds many different
possible notions of tameness agree. This indicates that
characterfree complex linear lagebraic groups
form
a good class of complex manifolds for the purpose of investigating
discrete subsets.

\begin{theorem}\label{main-1}
  Let $G$ be a character-free complex linear algebraic group
  and $\dim(G)\ge 2$. 
Let $D$ be a discrete subset of $G$.
Then the following statements are equivalent:
\begin{enumerate}
\item
  $D$ is (weakly) tame (in the sense of definition~\ref{def-tame}), i.e.,
  for every
exhaustion function $\rho:G\to\R^+$ and every map $\zeta:D\to\R^+$
there exists an automorphism $\phi$ of $G$ such that
$\rho(\phi(x))\ge \zeta(x)$ for all $x\in D$.
\item
  For every increasing positive sequence $r_k$ with $\lim r_k=\infty$
  there exists a biholomorphic self-map $\phi$ of $G$ such that
  \[
  \# \{ x\in \phi(D): \rho(x)\le r_k\}\le k  \ \forall k\in\N
  \]
\item
  $D$ is strongly tame (in the sense of definition~\ref{def-strong-tame}),
  i.e.,for every injective map $f:D\to D$ there exists
an automorphism $\phi$ of $G$ such that $\phi(x)=f(x)$ 
for all $x\in D$.
\item
  There is a biholomorphic self-map $\phi$ of the complex manifold
  $G$ and an algebraic subgroup $U$ of $G$
  such that $U\simeq(\C,+)$ (as a complex algebraic group)
  and $\phi(D)\subset U$.
\item
  There exists a biholomorphic self-map $\phi$  of $G$ and
  an algebraic subgroup $H$ of $G$ with $\dim(H)>0$ such that the quotient
  $Z=G/H$ is affine  and such that
  the natural projection map $\pi:G\to G/H=Z$ restricts
  to a proper map from $\phi(D)$ onto
  a discrete subset of $Z$.
\item
  There exists a biholomorphic self-map $\phi$  of $G$ and
  an algebraic subgroup $H$ of $G$ with $\dim(H)>0$
  such that the quotient
  $Z=G/H$ is affine and such that
  the natural projection map $\pi:G\to G/H=Z$ restricts
  to an injective map from $D$ onto
  a discrete subset of $Z$.
\item 
  There exists a biholomorphic self-map $\phi$ of $G$,
  a unipotent subgroup $U$ with $U\simeq(\C,+)$ and an open
embedding $i:G/U\hookrightarrow M$ into an affine variety $M$
such that $i(\pi(\phi(D)))$ is discrete in $M$ (where $\pi$ denotes
the natural projection $\pi:G\to G/U$.)
\item
  For every monotone increasing continuous function $h:[1,\infty[\to\R^+$
      with $\lim_{r\to\infty}\frac{h(r)}{\log r}=\infty$ there
      exists a biholomorphic self map $\phi$ of $G$ such that
      $N(r,\phi(D))< h(r)\ \forall r\ge 1$ where $N(r,D)$ is the counting
      function as defined in \S\ref{sect-nevanlinna}.
    \item
      For every discrete subset $D'$ of $G$ there exists a biholomorphic
      self-map $\phi$ of $G$ with $\phi(D)\subset D'$.
\end{enumerate}
\end{theorem}

\begin{proof}
  $(i)\implies(ii)$  If $D=\{a_k:k\in\N\}$, we may choose $\phi$ such that
  $\rho(\phi(a_k))>r_k$ for every $k$.

  $(i)\iff(iii)$ follows from Theorem~\ref{tame-eq}.

  $(iii)\implies(i)$: Obvious.

  $(i)\implies (iv)$: Due to Lemma~\ref{exist-unip} there
  exists an algebraic subgroup $U\simeq\C$. By Corollary~\ref{unbounded}
  every unbounded set (hence in particular $U$) contains a tame discrete
  subset $D'$ and due to Proposition~\ref{tame-eq} there is a holomorphic
  self-map $\phi$ of $G$ with $\phi(D)=D'\subset U$.

  $(iv)\implies(i)$: For dimension reasons we have $U\ne G$. Hence this
  implication follows from Proposition \ref{sub-uni}.
  
  $(i)\implies(vi)$: Theorem~\ref{torus-quot}.

  $(vi)\implies(v)$: Obvious.

  $(v)\implies(i)$: \cite{JW-TAME1},  Proposition~8.4.

  $(i)\implies (vii)$: Lemma~\ref{exist-unip}
  combined with Proposition~\ref{qa-quot}.

  $(vii)\implies(i)$: Corollary~\ref{overshear-tame}.

  $(i)\iff(viii)$: Proposition~\ref{equi-nevan}.
  
  $(i)\implies(ix)$:
  As an infinite discrete set, $D'$ is unbounded and therefore
  contains a tame discrete set $D''\subset D'$
  (Corollary~\ref{unbounded}). Now any two tame discrete subsets are
  equivalent (Corollary~\ref{cor-equiv}),
  therefore there is a biholomorphic self-map $\phi$ with
  $\phi(D)=D''\subset D'$.

  $(ix)\implies (i)$
  Consider a {\em tame} discrete subset $D'$. It follows directly from
  Definition~\ref{def-tame}
  that subsets of tame discrete sets are again tame. Hence
  $\phi(D)$ (and consequently $D$ itself) is tame.
  
\end{proof}
  
Remarks:

-- The implication $(i)\implies (iii)$ provides a positive
answer to question 7.7. of \cite{AU} for character-free
complex linear algebraic groups.
(Andrist and Ugolini asked for this implication for the class
of Stein manifold with volume density property - a class including
every character-free complex linear algebraix group.)

-- A discrete subset of a one-dimensional complex manifold
is never tame (\cite{JW-TAME1}); thus we need to require $G$ to have at least
dimension two.

-- Condition $(iii)$ implies in particular that every permutation
of a tame discrete subset extends to a biholomorphic self-map.

This assertion may be strengthened to the following result.
\begin{theorem}\label{extend-inj}
  Let $G$ be a character-free complex linear algebraic group
  and $\dim(G)\ge 2$. 

  Let $D$ and $D'$ be tame discrete subsets.

  Then every injective map $\zeta:D\to D'$ may be extended
  to a biholomorphic self-map of $G$.
\end{theorem}

Thus for example it follows from our results that every
injective map from $SL_n(\Z)$ to $SL_n(\Z[i])$ extends to a
biholomorphic self-map of $SL_n(\C)$.

\begin{corollary}\label{cor-equiv}
  Under the assumptions of the theorem any two tame discrete subsets
  $D$, $D'$ are ``equivalent'', i.e., there exists a biholomorphic
self-map of $G$ with $\phi(D)=D'$.
\end{corollary}

\begin{corollary}\label{cor-sym-diff}
  Let $G$ be a character-free complex linear algebraic group.
  Let $D,D'$ be infinite discrete subsets of $G$ such that
  the ``symmetric difference'' $(D\setminus D')\cup(D'\setminus D)$
  is finite.
  
  Then $D$ is tame if and only if $D'$ is tame.
\end{corollary}

\begin{proof}
  This follows immediately
  from $(v)$ or $(vii)$ of the theorem.
\end{proof}

\begin{corollary}
  Let $G$ be a character-free complex linear algebraic group.
  Let $D$ be a tame discrete subset. Let $A$ denote the group
  of all biholomorphic self-maps of $G$ fixing $D$ point-wise.

  Then $A$ acts transitively on $G\setminus D$.
\end{corollary}
\begin{proof}
  Let $p,q\in G\setminus D$. Due to
  Corollary~\ref{cor-sym-diff} both $D\cup\{p\}$ and $D\cup\{q\}$
  are tame.
  Define $\zeta:D\cup\{p\}\to D\cup\{q\}$ as $\zeta(x)=x$ for $x\in D$
  and $\zeta(p)=q$. Now Theorem~\ref{extend-inj} guarantees us that
  there exists a biholomorphic self-map $\phi$ of $G$ with $\phi|_D=id_D$
  (i.e.~$\phi\in A$) and $\phi(p)=q$.
\end{proof}

In the same way one deduces the following generalization:

\begin{corollary}
  Let $G$ be a character-free complex linear algebraic group and let
  $D$ be a tame discrete subset.

  Then the the automorphism group of the complex manifold $G\setminus D$
  acts $m$-transitively on $G\setminus D$ for every $m\in\N$,
  i.e., if $S$ and $S'$ are two subsets of $G\setminus D$, then
  there is a biholomorphic self-map of $G\setminus D$ mapping $S$
  to $S'$.
\end{corollary}

It should be noted that from the definition~\ref{def-tame} it is not obvious
that there exist any tame discrete subsets and indeed there are complex
manifolds without any tame discrete subset (\cite{JW-TAME1}).
Howeever, every character-free complex linear algebraic group of dimension
at least two does contain tame discrete subsets.
In fact, each of the characterizations $(iv)$--$(vii)$ of
theorem~\ref{main-1}
provides an easy way to construct tame
discrete subsets. In this way we obtain a confirmation of the following
known fact:
{\sl Every complex linear algebraic group of dimension at least two admits
  both tame and non-tame discrete subsets.}
(The existence of tame discrete subsets has been proved by Andrist and
Ugolini (\cite{AU}) for an arbitrary complex linear algebraic group $G$,
the existence of a non-tame discrete subset by the author for arbitrary
Stein manifolds (\cite{JW-LDS}).

  {\em History.} In \cite{JW-TAME1} the equivalence $(i)\iff(ii)$ was proved for
the special cases $G\simeq SL_n$ 
and $G\simeq(\C^n,+)$ 
and $(i)\iff (v)$ for $G\simeq SL_2$.

\subsection{Thrershold Sequences}
The technical notion of a ``Threshold sequence''
(which was first considered in \cite{JW-TAME1}) plays an important
role:
\begin{dfn}\label{def-threshold}
  Let $G$ be a complex manifold with exhaustion function $\rho$.
  A ``Threshold sequence'' for $G$ with respect to $\rho$ is
  an increasing sequence $R_k$ of positive real numbers such that
  the following property holds:
  {\em
    If $D\subset G$ is a discrete subset with
    \[
    \#\{z\in D:\rho(z)\le R_k\}\le k\quad \forall k\in\N
    \]
    then $D$ is tame.
  }
\end{dfn}

\begin{theorem}[=Theorem~\ref{threshold}]
  Let $G$ be a character-free complex linear algebraic group.
  
  Then $G$ admits a {\em threshold sequence} (for any exhaustion function
  $\rho$).
\end{theorem}

\subsection{Oka manifolds}

A complex manifold $X$ is called {Oka} iff for every $N\in\N$ and
every convex compact
subset $C\subset\C^N$ every holomorphic map with values in $X$ defined
on a neighbourhood of $C$ may be approximated by holomorphic maps
defined on all of $\C^N$.

(See \cite{Forstneric} for more information about Oka manifolds.)

Complex Lie groups are known to be Oka manifolds.
Now let $G$ be a complex linear algebraic group with
$Hom(G,\C^*)=\{1\}$. We show that certain manifolds
obtained from $G$ by certain modifications along tame discrete
subsets are still Oka:

\begin{theorem}
Let $G$ be a character-free complex linear algebraic group.

Let $D$ be a tame discrete subset.

Let $X=G\setminus D$ and let $Y$ denote the complex manifold obtained
from $G$ by blowing up $D$.
Then both $X$ and $Y$ are {\em Oka} manifolds.
\end{theorem}

\begin{proof}
  Theorem \ref{oka1} and Proposition~\ref{oka2}.
\end{proof}

{\em Remark.} For the class of character-free
complex linear algebraic groups this answers a question raised by
Andrist and Ugolini (\cite{AU}, Question 7.8) for
Stein manifolds with volume density property.

In contrast, every Stein manifold $X$
(in particular, every complex linear algebraic group)
admits a discrete subset $D$ such
that $X\setminus D$ is not Oka. In (\cite{JW-LDS}) it is proved:
Every Stein manifold $X$ admnits a discrete subset $D$ such that
every holomorphic map from $\C^n$ (with $n=\dim(X)$) to $X\setminus D$
is degenerate,
a property which implies that $X\setminus D$ can not be an Oka manifold.

\subsection{Discrete subgroups}

We conjecture that every discrete subgroup of a complex linear algebraic
group is a tame discrete set.
In this direction we have some partial results.

\begin{theorem}
Let $G$ be a complex linear algebraic group and let $\Gamma$  be a discrete
subgroup which contains a unipotent element other than the neutral
element. 

Assume that furthermore
at least one of the following conditions is fulfilled:
\begin{enumerate}
\item
$G\simeq SL_2(\C)$ or $G\simeq PSL_2(\C)$.
\item
$G$ is defined over $\Q$ or a totally imaginary quadratic number
field $K$ and $\Gamma$ is contained in $G(\O)$ where $\O$ denotes
the ring of algebraic integers of $\Q$ resp.~$K$.
\end{enumerate}

Then $\Gamma$ is a tame discrete subset of the complex manifold $G$.
\end{theorem}

The theorem implies in particular that $SL_n(\Z)$ and 
$SL_n(\Z[i])$ are tame discrete subsets
of $SL_n(\C)$. 

In the context of this theorem
it is an interesting fact that, if $\Gamma$ is a
discrete subgroup in a complex linear algebraic group $G$ such that
$G/\Gamma$ has finite volume with respect to the Haar measure of $G$,
but is not compact, then $\Gamma\setminus\{e\}$
contains a unipotent element.
(Theorem of Kazhdan-Margulis, \cite{KM}.)

Thus, if $G=SL_2(\C)$ and $\Gamma$ is a discrete subgroup,
such that $G/\Gamma$ has finite volume, but is not compact, then
$\Gamma$ is a tame discrete subset.

If $\Gamma$ is the fundamental group of a hyperbolic punctured Riemann
surface, then $\Gamma$ embedded into $PSL_2(\R)$ by its action of
deck transformations on the upper half plane becomes a tame discrete
subset of $PSL_2(\C)$. 

Other examples are provided by fundamental groups of hyperbolic knot
complements. Except for certain special cases, the complement
$M=S\setminus C$
of a knot (or link)
$C$ in the $3$-sphere $S$ admits a complete Riemannian metric of
negative constant sectional curvature, i.e., $M$ is a ``real hyperbolic
manifold''. This implies that the universal covering $\tilde M$ of $M$
is isomorphic to $H=G/K$ with $G=SO(3,\C)$ and $K=SO(3,\R)$.
Moreover, $G$ is the group of isometries of the ``hyperbolic $3$-space''
$G/K$. Let $\Gamma$ be the fundamental group of $M$. Then $\Gamma$
acts by ``deck transformations'' on its universal covering
$\tilde M =G/K$. Thus $\Gamma$ admits a natural embedding into
$G=SO(3,\C)$ which is the group of isometries of $\tilde M$.
In this way $\Gamma=\pi_1(M)$ is naturally embedded into the
complex linear algebraic group $G=SO(3,\C)$ as a discrete
subgroup. By our results it follows that this is a {\em tame}
discrete subset. (see \S\ref{sect-hyper}.)

\subsection{Union of tame sets}

\begin{theorem}[=Theorem~\ref{union-2-tame}]
Let $D$ be a discrete subset of a linear algebraic group $G$ with
$\dim(G)\ge 2$.

Then $D$ can be realized as the union of two tame subsets
of the complex manifold $G$.
\end{theorem}

\section{Relation to Nevanlinna Theory}\label{sect-nevanlinna}
Here we want to explain the connection between tameness and the usual
notions of Nevanlinna theory.
(See e.g.~\cite{NW} for more information about Nevanlinna theory.)

Given a divisor $D=\sum_k n_k\{p_k\}$ in $\C$,
  in Nevanlinna theory
  one introduces the ``counting function'' $N(r,D)$
  which may be defined as
  \[
  N(r,D)=
  \int_1^r n(t,D) \frac{dt}{t}\quad \text{ with }
  n(t,D)=\sum_{k:|p_k|\le t}n_k
  \]

  In this spririt we define for a
  discrete subset $D$ in a complex manifold
  $X$ with exhaustion function $\rho$
    \[
    N(r,D)=\int_1^r n(t,D) \frac{dt}{t}\quad \text{ with }
  n(t,D)=\#\{x\in D:\rho(x)\le t\}
  \]

  We recall the notation $\log^+(x)\stackrel{def}=\max\{x,0\}$.

  \begin{remark}
    We may replace a given exhaustion function $\rho$ by
    $\tilde\rho(p)=\max\{1,\rho(x)\}$ and therefore without
    loss of generality assume that $\rho\ge 1$. Under this assumption
    we have
    \[
    N(r,D)=\sum_{p\in D}\log^+\frac{r}{\rho(p)}
    \]
  \end{remark}

  The goal of this section is to establish the following result.
  
  \begin{proposition}\label{equi-nevan}
    Let $X$ be a complex manifold.
    
    A discrete subset $D$ in $X$ is tame if and only if
    for every real continuous increasing function $h:[1,\infty[\to\R^+$
        with $\lim_{r\to\infty}\frac{h(r)}{\log r}=+\infty$
        there exists a biholomorphic self-map $\phi$ of $X$ such that
        $N(r,\phi(D))\le h(r)\ \forall r\ge 1$.
  \end{proposition}

  \begin{proof}
    This is an immediate consequence of Corollary~\ref{R2h} and
    Proposition~\ref{h2R} below.
  \end{proof}

  \begin{lemma}
    Let $h:[1,\infty)\to\R^+$ be a continuous increasing function
      with $\lim_{r\to\infty} h(r)/\log(r)=+\infty$.
      Then there exists a positive real number $R>1$
      such that
      \[
      \log^+ \frac{r}{R}< h(r)\ \forall r\ge 1
      \]
  \end{lemma}
  \begin{proof}
    Since $\lim_{r\to\infty} h(r)/\log(r)=+\infty$, we have
    \[
    \lim_{r\to\infty}\left(  \log r -h(r) \right) =-\infty
    \]and therefore there exists an upper bound
    $B$ for $\log r-h(r)$ on $[1,+\infty)$.
      Now it suffices to choose $R$ such that $\log R>B$.
  \end{proof}
  \begin{corollary}\label{cor-3-2}
  Let $h:[1,\infty)\to\R^+$ be a continuous increasing function
      with $\lim_{r\to\infty} h(r)/\log(r)=+\infty$.
      Then there exists a sequence of  positive real numbers $R_k$
      with $\lim_{k\to\infty} R_k =+\infty$ such that
      \[
      \sum_k \log^+\frac{r}{R_k}< h(r)\ \forall r\ge 1
      \]
  \end{corollary}
  \begin{proof}
    For each $k\in\N$ we apply the lemma to $h_k(r)=2^{-k}h(r)$
    in order to obtain a real number $R_k$ such that
    \[
    \log^+\frac{r}{R_k}< h_k(r)= 2^{-k}h(r)\ \forall r\ge 1
    \]
    Since this inequality is only strengthened by increasing $R_k$
    there is no problem in choosing the constants $R_k$ such that
    $\lim_{k\to\infty}R_k=+\infty$.
    
    Summing up over $k$ yields the desired inequality. (We observe that
    due to $\lim_{k\to\infty}R_k=\infty$ the sum $\sum_k \log^+(r/R_k)$
    is actually locally finite. Hence there are no convergence problems.)

  \end{proof}

  \begin{corollary}\label{h2R}
    Let $h:[1,\infty[\to\R^+$ be an increasing continuous function
        with $\lim_{r\to\infty}h(r)=+\infty$.

        Then there exists an increasing
    sequence of positive real numbers $R_k\ge 1$
    with $\lim R_k=\infty$ such that for every discrete subset $D\subset X$
    the property
    \begin{equation}\label{eqstar}
          \#\{x\in D: \rho(x)\le R_k\}\le k\ \forall k
    \end{equation}
          implies $N(r,D)\le h(r)\ \forall r\ge 1$.
  \end{corollary}
  \begin{proof}
    Due to Corollary~\ref{cor-3-2} there exists such a sequence $R_k$
    with
    \[
      \sum_k \log^+\frac{r}{R_k}< h(r)\ \forall r\ge 1
      \]
      On the other hand, \eqref{eqstar} implies
      \[
      N(r,D)\le \sum_k \log^+\frac{r}{R_k}
      \]
      Hence the statement.
  \end{proof}
  \begin{proposition}\label{R2h}
    Let $R_k$ be an increasing
    sequence of positive real numbers $R_k\ge 1$
    with $\lim R_k=\infty$.
    Then there exists an increasing continuous function $h:[1,\infty[\to\R^+$
        with $\lim_{r\to\infty}h(r)/\log r=\infty$
        such that the following implication holds for every discrete set $D$:

        {\em If $N(r,D)\le h(r)$ for all $r\ge 1$, then
          \[
          \#\{x\in D: \rho(x)\le R_k\}\le k\ \forall k
          \]
        }
  \end{proposition}
  \begin{proof}
    Let $R_k$ be a given real positive sequence with $\lim R_k=+\infty$.
    We choose $c_1>1$ such
    \[
    \frac{2c_1}{c_1+\log R_1}>\frac12.
    \]
    (This is possible, since $\lim_{c\to\infty} 2c/(c+\log R_1)=2$.)
    Then we recursively choose $c_k$ for every $k\in\N\setminus\{1\}$
    such that
  \begin{equation}\label{ck}
  1 < c_k,\quad\quad  c_{k-1}<c_ke^{c_{k-1}}R_{k-1}<e^{c_k}R_k,
  \quad\quad \frac{c_k(k+1)}{c_k+\log R_k}>\frac k2
  \end{equation}
  Next we
  observe that $e^{c_k}R_k$ is strictly increasing and
  define a continuous function $h:[1,\infty(\to\R^+$ by
    stipulating
        the following conditions:
    \begin{enumerate}
    \item
      $h(e^{c_k}R_k)=c_k(k+1)$ and
    \item
      For every $k\in\N$ there are constants $a_k,b_k$ such that
      $h(t)=a_k+b_k\log t$ on the interval
      $[e^{c_k}R_k,e^{c_{k+1}}R_{k+1}]$.
    \end{enumerate}
    Define $p_k=e^{c_k}R_k$.
    By construction we have $h(p_k)=c_k(k+1)$, implying
    $h(p_k) < h(p_{k+1})\ \forall k$.
    It follows that $h(t)=a_k+b_k\log t$ with $b_k>0$
    on $[p_k,p_{k+1}]$. Thus $h$ 
    is monotonously increasing.
    We also observe that
    \[
    \frac{h(t)}{\log t}=\frac{a_k+b_k\log t}{\log t}
    =\frac{a_k}{\log t}+b_k
    \]
    is monotone on each interval $[p_k,p_{k+1}]$.
    Since
    \[
    \frac{h(p_k)}{\log p_k}=\frac{c_k(k+1)}{\log p_k}
    =\frac{c_k(k+1)}{c_k+\log R_k}> \frac{k}{2},
    \]
    it follows that $\lim_{t\to\infty} h(t)/\log t=+\infty$.

    Now let us assume that there exists a number $k$
    with
    \begin{equation}\label{notame}
    \#\{x\in D: \rho(x)\le R_k\}> k
    \end{equation}
    This implies that $n(r,D)\ge k+1$ for all $r\ge R_k$.
    Hence
    \begin{multline*}
    N(p_k,D)=N(e^{c_k}R_k,D)\ge \int_{R_k}^{e^{c_k}R_k}n(t,D)dt/t\ge\\
    \int_{R_k}^{e^{c_k}R_k}(k+1)dt/t= c_k(k+1)=h(p_k).
    \end{multline*}
    Thus the existence of a natural number $k$ with \eqref{notame}
    implies that there exists a number $r\ge 1$ with $N(r,D)\ge h(r)$.
  \end{proof}
    
  \section{Tools}
Here we summarize some mostly well-known material which we will need.

\begin{fact}
Let $G$ be an algebraic group with algebraic subgroups $H$ and $I$.

Then the following conditions are equivalent:
\begin{enumerate}
\item
The $H$-action (by left multiplication) on $G/I$ is trivial.
\item
$H\subset gIg^{-1}$ for every $g\in G$.
\item
$gHg^{-1}\subset I$ for every $g\in G$.
\end{enumerate}
\end{fact}

\begin{fact}
Let $G$ be a linear-algebraic group and let $H$ be an
algebraic subgroup.

\begin{enumerate}
\item
If $H$ is normal or reductive, then $G/H$ is an affine variety.
(\cite{MM}).
\item
If every homomorphism of algebraic groups from $H$ to $\C^*$
is trivial (e.g., if $H$ is unipotent), then
$G/H$ is a quasi-affine variety.
\end{enumerate}
\end{fact}

\begin{fact}
  Let $D$ be a discrete subset of a Stein manifold $X$, let $G$
  be a complex Lie group and let $\alpha:D\to G$ be a map.

  Then there exists a holomorphic map $F:X\ss to G$ with
  $F(x)=\alpha(x)$ for all $x\in D$.
\end{fact}

\begin{lemma}[ nil cone]\label{nil-cone}
  Let $G$ be a character-free complex linear algebraic group.
  Let $N$ denote the ``nil cone'' defined as the set of all
  elements $v$ of the Lie algebra $Lie(G)$ for which $ad(v)$
  is nilpotent, i.e., 
\[
N=\{v\in \Lie(G): \exists n : ad(v)^n=0\}.
\]

  Then $N$ is an irreducible analytic subset of $\Lie(G)$
  which generates $\Lie(G)$ as a complex vector space.
\end{lemma}

\begin{proof}
  An element $v\in\Lie(G)$ is nilpotent if and only if
  the corresponding one-parameter group
  \[
  \{\exp (vt):t\in\C\}
  \]
  is a unipotent subgroup of $G$.
  Let $B$ be a Borel subgroup of $G$ und let $U$ denote the unipotent
  radical of $B$. Every unipotent subgroup of $G$ is conjugate to
  a subgroup of $U$. Therefore $N$ is the set of all elements of $\Lie(G)$
  which can be written as $c=Ad(g)(v_0)$ with $g\in G$, $v_0\in\Lie(U)$.
  This proves irreducibility of $N$.

  The nil cone $N$ is evidently invariant under the adjoint action.
  It follows that the vector space $V$ spanned by $N$ is likewise
  invariant under the adjoint action and therefore an ideal of
  the Lie algebra $\Lie(G)$.
  However, structure theory of linear algebraic groups implies  that
  the only normal Lie subgroup of $G$ containing $U$ is $G$ itself.
  Hence $V=\Lie(G)$.
\end{proof}

\begin{corollary}
\label{exist-unip}
Let $G$ be a character-free complexlinear algebraic group.
Then $G$ contains an algebraic subgroup $U\simeq(\C,+)$.
\end{corollary}
\begin{proof}
  Define $U=\{\exp(vt):t\in \C\}$ for some element $v\in N\setminus\{0\}$.
\end{proof}

\begin{lemma}\label{unbounded-orbit}
Let $G$ be a complex linear-algebraic group, and let $I$ be an algebraic 
subgroup such that $G/I$ is quasi-affine. Let $i:G/I\to Y$ be an
open dense embedding into an affine variety $Y$.

Let $H$ be an algebraic subgroup of $G$ and let $M_0=Hp$ be a
positive-dimensional orbit in  $G/I$.

Then $M=i(M_0)$ can not be relatively compact in $Y$.
\end{lemma}

\begin{proof}
Assume the contrary, i.e., assume that $M$ is relatively compact in $Y$.
Since $Y$ is affin, points on $Y$ may be separated by regular functions.
Hence there must be a regular function $f$ which is not constant on $M$.
Since $M$ is assumed to be relatively compact in $Y$,
every continuous function on $Y$
is bounded on $M$.

Now consider the exponential map $\exp:\Lie(H)\to H$.
Via $F:v\mapsto f(\exp(v)\cdot p)$ we obtain a holomorphic function
on the complex vector space $\Lie(H)$.
Note that $F$ is bounded, because $f$ is bounded. Hence $F$ must be
constant due to the theorem of Liouville. Since $\exp$ is locally 
biholomorphic near the origin in $\Lie(H)$, it follows that $f$ must be
constant on $M$ which yields a contradiction. Thus $M$ can not be
relatively compact in $Y$.
\end{proof}

\begin{lemma}\label{4.2}
  Let $G$ be a character-free complex linear algebraic group,
  $T$ a maximal torus, and $U$ a $T$-invariant maximal unipotent subgroup
(i.e., $U$ is maximal unipotent and $gUg^{-1}=U$ for all $g\in T$).
Then there exists an dense open embedding of $G/U$ into an affine variety $V$,
and a point $p\in G/U$ such that
\begin{enumerate}
\item
The isotropy
group $T_p=\{x\in T: x(p)=p\}$ is trivial and
\item
 the $T$-orbit $T(p)$
is closed in $V$.
\end{enumerate}
\end{lemma}

\begin{proof}
By taking the quotient by the unipotent radical we may assume that
$G$ is semisimple. 
We choose a faithful irreducible representation $\rho$ of $G$. 
Thus we have an injective morphism of algebraic groups 
$\rho:G\to SL(n,\C)$ for some $n$. We may assume that $\rho(T)$ is
contained in the diagonal matrices, because $\rho(T)$ is contained in
a maximal torus in $SL(n,\C)$, hence up to conjugation in the set $D$
of diagonal matrices of determinant $1$.
We may also assume that $\rho(U)$ is contained in the set $U^+$ of upper
triangular matrices.
We consider $\C^n\setminus\{(0,\ldots,0)=SL_n(\C)/H$,
with
\[
H= \left\{
\left(
\begin{array}{c|c}
  1 & v \\
  \hline
  0 & A \\
\end{array}
\right)
:
A\in SL_{n-1}(\C), v\in\C^{n-1}
\right\}
\]
We note that
$U^+\subset H$.
Define 
\[
\Omega=\{z\in\C^n: z_i\ne 0\ \forall i\}.
\]
We consider
\[
G/U\longrightarrow SL_n(\C)/H \simeq\C^n\setminus\{(0,\ldots,0)\}
\]
with the base point $eU$ being mapped to $e=(1,0,\ldots,0)$.
Since the representation $\rho$ was assumed to be irreducible,
the $G$-orbit through $e$ intersects $\Omega$.
Hence there is an element $g_0\in G$ such that $g_0U$ is mapped
to a point $q$ in $\Omega$.  Observe that the $D$-orbit through $q$
is closed in $\C^n$, because it equals
\[
\{(z_i,\ldots,z_n):\Pi z_i=c\}
\]
for some $c\in\C^*$.
Recall that $\rho$ embedds $T$ as a Zariski closed subgroup into $D$.
It follows that the $T$-orbit through $q$ is closed as well.
Moreover the $T$.isotropy group at $q$ is trivial, because the isotropy
group of $D$ at $q$ is trivial.
Now let $\zeta:G/U\to W$ be an open embedding into an affine variety
$W$. (Note that $G/U$ is quasi-affine, because $U$ is unipotent.)
The coordinates $z_i$ on $\C^n$ induce regular functions $\xi_i$ on $G/U$
via
\[
gU \mapsto z_i(\rho(g)(p)).
\]
We define 
\[
F:G/U\to W\times\C^n
\]
 as 
\[
F(gU)=(\zeta(gU),\xi_1(gU),\ldots,\xi_n(gU))
\]
and take $V$ to be the closure of $F(G/U)$ in $W\times\C^n$.
Define $p=g_0U$. Note that the $T$-orbit through $p$ is closed and
with trivial isotropy. Hence the statement. 
\end{proof}

\begin{proposition}\label{jw1-maxtorus}
Let $S$ be  a complex semisimple Lie group with exhaustion function 
$\rho:S\to\R^+$. Then there exists
an increasing sequence $R_n$
of positive real numbers such that the following holds:

{\em ``If $D$ is a discrete subset of $S$ with
\[
\#\{x\in D:\rho(x)\le R_n\}\le n\ \forall n\in\N
\]
then there exists a maximal torus $T$ and a biholomorphic self map $\phi$
of $S$ such that the natural projection map $\pi:S\to S/T$
restricts to an injective proper map from $\phi(D)$ to $S/T$.''}
\end{proposition}
\begin{proof}
\cite{JW-TAME1}
\end{proof}

\begin{proposition}\label{proj-tame}
  Let $G$ be a complex Lie group, $H$ a closed complex Lie subgroup
  $D$ a discrete subset of $G$, and let $\pi:G\to G/H$ be the
  projection map.

  Assume that $G/H$ is Stein, and that the projection map $\pi$
  restricts to a finite map from $D$ to a discrete subset of $G/H$.

  Then $D$ is tame.
\end{proposition}
\begin{proof}
  See \cite{JW-TAME1}
\end{proof}

An easy modification of the proof given in \cite{JW-TAME1}
yields:
\begin{proposition}\label{4.6}
  Let $G$ be a linear-algebraic group with an algebraic subgroup
  $H$, such that the quotient variety $G/H$ admits an open dense
  embedding $i:G/H\to M$ into an affine variety $M$.
  Let $D$ be a discrete subset.

  If the natural projection morphism $\pi:G\to G/H$ composed with the
  embedding $i:G/H\to M$ restricts to a finite map from $D$ to a
  discrete subset of $M$, then $D$ is tame.
\end{proposition}

\begin{proposition}
  Let $G$ be a character-free complex linear algebraic group.
  
Let $T$ be a connected reductive commutative subgroup of $G$.

Then every discrete subset of $T$ is a {\em tame} discrete subset
of $G$.
\end{proposition}

\begin{proof}
Lemma~\ref{4.2} combined with Proposition~\ref{4.6}.
\end{proof}

\section{$Aut(X)$ is a Baire space}
\begin{lemma}\label{aut-baire}
The automorphism group of a complex manifold $X$, endowed with
the compact-open-topology, is a Baire space.
\end{lemma}

\begin{proof}
Let $d_0$ be a complete metric on $X$. Then $d_1=\min\{1,d_0\}$ is a bounded
complete metric. Let $\Omega_k$ be a countable family of
relatively compact open subsets of $X$ which cover all of $X$. 
We define a metric on the automorphism group $G$ of $X$ as follows
\[
d(f,g)=\sum_k \frac{1}{2^k}
\sup_{x\in\Omega_k}
\left( \max\{ d_1((f(x),g(x)), d_1(f^{-1}(x),g^{-1}(x)) \} \right)
\]
This turns $G$ into a complete metric space.
It is well-known that
complete metric spaces are Baire spaces.
\end{proof}

\begin{lemma}
$Aut(X)$ has a countable basis of topology because $X$ does.
\end{lemma}

\begin{proof}
Let $(U_k)_{k\in\N}$ be a countable basis of the topology of $X$.
Assume that all the $U_k$ are relatively compact.
Let $K$ be compact and $W$ open in $X$. We claim that for any
$f:X\to X$ with $f(K)\subset W$ we can find a map
\[
\N\ni i\mapsto (j_{i,1},\ldots,j_{i,s_i}) \ \ s_i\in\N_0
\]
such that
\begin{enumerate}
\item
$f(\overline{U_i})\subset\cup_{1\le k\le s_i}U_{j_{i,k}}\ \forall i$
(with $\cup_{k\in\{\}}U_k=X$.)
\item
$K\subset\cup_{i:s_I>0}\overline{U_i}$
\item
$U_{j_{i,k}}\ \forall i,1\le k\le s_i$.
\end{enumerate}
\end{proof}

\begin{lemma}\label{6-3}
Let $X$ be a complex manifold on which its automorphism group 
$G=\Aut(X)$ acts transitively.
Let $p\in X$ and let $Z\subset X$ be a nowhere dense closed
analytic subset.

Then $B=\{g\in G: g(p)\not\in Z\}$ is open and dense in $G$.
\end{lemma}

\begin{proof}
Openness of $B$ is obvious.
  
  Let $S$ be a countable dense subset of $G$. If $B$ is not dense,
  then $A=G\setminus B$ has non-empty
interior, there is an open subset $W\subset G$ with $Wp\subset Z$.
But then $SW=G$, implying $Gp\subset\cup_{s\in S}s\cdot Z$ which contradicts
the assumption of transitivity (because $S$ is countable and $Z$ is
a nowhere dense analytic subset).
\end{proof}

\begin{corollary}\label{prep-baire}
  Let $X$ be a complex manifold on which its automorphism group $Aut(X)$
  acts transitively. Let $Z$ be a nowhere dense closed analytic subset of $X$
  and let $D$ be a countable subset of $X$.

  Then there exists an automorphism $g$  of $X$ with
  $g(p)\not\in Z\ \forall p\in D$.
\end{corollary}

\begin{proof}
  Due to Lemma~\ref{6-3}, the set $U_p=\{g\in Aut(X):g(p)\not\in Z\}$ is
  open and dense in $Aut(X)$ for every $p\in D$. Now $Aut(X)$ is Baire
  (Lemma~\ref{aut-baire}). Hence the intersection of countably many
  dense open subsets of $Aut(X)$ is not empty. This implies the claim.
\end{proof}

\section{Projections}

A key method in this theory is to exploit the interplay between
two different projection maps.
For this method the proposition below is key technical result.

\begin{proposition}\label{pingpong}
Let $G$ be a linear-algebraic group with discrete subsets $D$, $D'$ and
algebraic subgroups $H$ and $I$.

Assume that $X= G/H$ admits an open embedding into an affine variety 
$X'$
and that $Y=G/I$ admits an open embedding into an affine variety $Y'$.

Furthermore we assume that the natural projection $\pi:G\to X=G/H$ 
restricts to a proper map
from $D$ and $D'$ to $X'$, i.e.,  the fibers are finite and  the 
images $\pi(D)$ and $\pi(D')$
are discrete in $X'$.

Let $\tau:G\to Y=G/I$ denote the natural projection.

Let $\zeta:D\to D'$ be a bijective map.

\begin{enumerate}
\item
If the $H$-left action on $G/I$ is non-trivial, then there exists
a holomorphic automorphism $\phi$ of the complex manifold $G$
such that $\tau\circ\phi$  maps
$D$ to a discrete subset of $Y'$ with finite fibers.
Furthermore $\phi$ can be chosen such that
$\tau\circ\phi$ is injective if $H$ acts freely on $G/I$ or
$\pi|_D$ is injective.
\item
Assume that the $H$-action on $G/I$ has a dense orbit
and that the projection map to $G/H$ is injective on $D\cup D'$.
Then there exists a holomorphic automorphism $\phi$ of the 
complex manifold $G$
such that $\tau\circ\phi$ is injective on $D'$,
$\tau\circ\phi=\tau\circ\phi\circ\zeta$ on $D$
and
such that $\tau(\phi(D'))=\tau(\phi(\zeta(D)))$ is discrete
 in $Y'$.
\end{enumerate}
\end{proposition}

\begin{proof}

Let $D=\{a_k:k\in\N\}$ and $D'=\{b_k:k\in\N\}$ with $\zeta(a_k)=b_k$.

We assume that the $H$-action on $G/I$ is non-trivial.
(We may do this, because the statement is empty if $H$ acts
trivially on $G/I$.)

Let $E_0$ denote the set of all points in $Y$ where the dimension of
the $H$-orbit is less than the generic $H$-orbit dimension. 
Note that $E_0$ is a nowhere dense closed  analytic subset of $Y$.
Let
$E=\tau^{-1}(E_0)$. $E$ is a closed analytic set, $E\ne G$
 and $D\cup D'$ is
countable. Therefore there exists an element $g\in G$ such that
$ga_k,gb_k\not\in E\ \forall k$, i.e., we may wlog assume that
$(D\cup D')\cap E=\{\}$.

$(i)$ Since $D\cap E=\{\}$, the $H$-orbit through $\tau(a_k)$ is
non-trivial for every $k$. Now every non-trivial $H$-orbit is unbounded
(i.e.~has non-compact closure)
in $Y'$
(lemma~\ref{unbounded-orbit}). 
Hence we may choose a sequence $h_k\in H$ such that
$h_k\cdot\tau(a_k)=\tau(h_k\cdot a_k)$ is a discrete sequence in $Y'$.
Since $\pi$ restricts to a map with finite fibers on $D$,
we may choose this sequence in such a way that furthermore
$h_k=h_m$ whenever $\pi(a_k)=\pi(a_m)$.  If $\pi|_D:D\to X$ is injective,
we chose $h_k$ in such a way that $k\mapsto h_k\cdot\tau(a_k)$ is
injective, too.
Observe that $\pi(a_k)\mapsto h_k$ is a well-defined map from
the discrete subset $\pi(D)\subset X'$ to $H$. Because $H$
is a complex Lie group and $X'$ is affine (hence Stein), we can
find a holomorphic function $f:X'\to H$ such that
$f(\pi(a_k))=h_k$ for every $k\in\N$. We define a biholomorphic
self-map of $G$ by
\[
\phi: x\mapsto x\cdot f(\pi(x))
\]
This proves statement  $(i)$
$(ii)$
We assume that $H$ has a dense orbit in $G/I$. Due to algebraicity
this dense orbit is open.
For every $k$ both $\tau(a_k)$ and $\tau(b_k)$ are contained
in this (unique) open $H$-orbit in $Y$. Hence we may choose 
$h_k,\tilde h_k\in H$
such that $\tau(h_ka_k)=\tau(\tilde h_kb_k)$ and such that this sequence
$\tau(h_ka_k)$ is discrete in $Y'$.

Next we observe that $\alpha:g\mapsto g^{-1}$ is a biholomorphic self-map of 
$G$ inducing an isomorphism between the left quotient $H\backslash G$
and the right quotient $G/H$. Therefore we may assume that the
natural projection $\tilde\pi:G\to H\backslash G=Z$ maps $D$ injectively
onto a discrete subset of $Z'$ where $Z'$ is an affine variety 
containing $H\backslash G$ as a Zariski open subvariety.

In both cases it suffices now to choose a holomorphic map $f$ from
$Z$ to $H$ such that $f(\pi(a_k))=h_k$ and $f(\pi(b_k))=\tilde h_k$.

Then the automorphism defined by
\[
\phi: g \mapsto f(Hg)\cdot g 
\]
does the job.
\end{proof}

\begin{theorem}\label{torus-quot}
Let $G$ be a character-free complex linear algebraic group.
Let $D$ be a tame discrete subset.
Let $T$ denote the maximal torus of $G$.

Then there exists a biholomorphic self-map $\phi$  of $G$
such that the projection map $\pi:G\to G/T$ restricts to an
injective proper map on $\phi(D)$.
\end{theorem}

\begin{proof}
Let $R$ denote the unipotent radical of $G$ and let $S$ be a maximal
semisimple algebraic subgroup of $G$.
If $R=\{e\}$, then $G$ is semisimple and we may invoke
proposition~\ref{jw1-maxtorus}.
If $S=\{e\}$, then $T=\{e\}$ and the statement becomes trivial.
Thus we may assume that both $R$ and $S$ 
are positive-dimensional.

We observe that group multiplication
yields a bijective map $R\times S\to G$ given by $(x,y)\mapsto x\cdot y$.
Let $a_k$, $b_k$ be sequences in $R$ resp.~$S$ such that
\[
D=\{a_k\cdot b_k: k\in\N\}.
\]
Since all maximal tori and all maximal semisimple subgroups are
conjugate, we may assume that $T\subset S$.
We fix an exhaustion function $\rho_S$ on $S$.
As a variety, $R$ is isomorphic to $\C^n$ and we may use the
euclidean norm to define an exhaustion function $\rho_R$
on $R$ in such a way that there is a holomorphic action of a compact
Lie group $K$ (namely $SU(n,\C)$)
on the complex manifold $R$ which leaves $\rho_R$
invariant.
Now we define an exhaustion function $\rho$ on $G$ by
imposing 
\[
\rho(xy)=\rho_R(x)+\rho_S(y)\ \forall x\in R,y\in S.
\]
Let $R_k$ be a threshold sequence for $S$ with respect to the
exhaustion function $\rho_S$ (cf. \cite{JW-TAME1}, prop.~9.6). 
Since we assume $D$ to be tame,
we may assume that 
\[
\#\{x\in D:\rho(x)\le 2R_k+1\}\le k\ \forall \in\N.
\]
Applying a small translation if necessary, we may assume
that $D\cap S$ is empty while
\begin{equation}\label{eqrk}
\#\{x\in D:\rho(x)\le 2R_k\}\le k\ \forall \in\N.
\end{equation}
holds.

Now we
define
\begin{align*}
D' &=\{a_k\cdot b_k: k\in\N, \rho_R(a_k)>\rho_S(b_k)\} \\
D''&=\{a_k\cdot b_k: k\in\N, \rho_R(a_k)\le \rho_S(b_k)\} \\
\end{align*}
This definitions imply that
\[
\{a_k: k\in\N, (a_k,b_k)\in D'\}
\]
is discrete in $R$ and
\[
\{b_k: k\in\N, (a_k,b_k)\in D''\}
\]
is discrete in $S$.

For $(a_k\cdot b_k)\in D''$ we have $\rho_R(a_k)\le \rho_S(b_k)$
implying
\[
\rho_S(b_k)\ge \frac 12\rho(a_k\cdot b_k)
=\frac 12 \left(
\rho_R(a_k)+\rho_S(b_k)
\right)
\]
which in combination with  \eqref{eqrk} implies
\[
\#\{a_k: k\in\N, (a_k,b_k)\in D'',\rho_S(b_k)\le R_k\}\le k\ \forall k\in\N.
\]
It follows that
\[
\{a_k: k\in\N, (a_k,b_k)\in D''\}
\]
is a tame discrete subset of $S$.
Hence, due to proposition~\ref{jw1-maxtorus}, 
the projection map $R\times S\to S/T$ restricts
to an injective proper map from $D''$ to $S/T$
(after applying a suitable biholomorphic self-map of $S$).

Now $D''$ maps  properly to $S/T$, while $D'$ maps properly to $R$.

Put together, these facts imply that $D=D'\cup D''$ maps
properly to $G/T\simeq (R\times S)/T$.

We observe that the intersection of a reductive and a unipotent subgroup
is necessarily trivial. Therefore $T$ acts freely on $G/R$ and $R$ acts
freely on $G/T$. Using proposition~\ref{pingpong} (with
$T$ as $H$ and $R$ as $I$) we obtain that - after applying a suitable
biholomorphic self-map - the projection of $D$ to $G/R$ is proper and
injective. Now another invokation of proposition~\ref{pingpong}
(this time with $R$ as $H$ and $T$ as $I$) yields the statement.
\end{proof}

  The above prove provides some additional information.
  In fact, we have shown:
{\em Let $G$ be a complex linear algebraic group with maximal semsimple
    subgroup $S$ such that every morphism of algebraic groups
    from $G$ to $\C^*$ is trivial.
    $\dim(S)>0$. If $R_k$ is a threshold sequence for $S$,
then $2R_k+1$ is a threshold sequence for $G$.}

As a consequence we obtain:

\begin{proposition}\label{threshold}
  Let $G$ be a character-free complex linear algebraic group
  with $\dim(G)\ge 2$.
  
  Then there exists a ``threshold sequence'' $R_k$.
\end{proposition}

\begin{proof}
  If $G$ is not solvable, then $¿dim(S)>0$ and
  the statement  follows from the above considerations.
  If $G$ is solvable, then the triviality of $\Hom(G,\C^*)$ implies
  that $G$ is unipotent. In particular, $G\simeq\C^n$ as an algebraic
  variety. Then the existence of a threshold sequence follows from
  \cite{JW}.
\end{proof}

\begin{corollary}\label{unbounded}
  Let $G$ be a complex linear algebraic group, $\dim(G)\ge 2$ such that
  every morphism of algebraic groups to $\C^*$ is trivial.

  Then every unbounded subset $A$ of $G$ contains a tame discrete
  subset of $G$.
\end{corollary}

\begin{proof}
  Let $\rho$ be an exhaustion function on $G$.
  By Proposition~\ref{threshold} there exists a threshold sequence $R_k$ for $G$
  with respect to $\rho$. 
  If $A$ subset $A$ of $G$ is unbounded, we may choose a sequence of points
  $p_k$ of $A$ with $\rho(p_k)>R_k$.
  Now $D=\{p_k:k\in\N\}$ is a tame discrete subset of $G$ which is
  contained in $A$.
\end{proof}

\begin{proposition}\label{qa-quot}
Let $G$ be a character-free complex linear algebraic group.
Let $D$ be a tame discrete subset.
Let $H$ be an algebraic subgroup of $G$ such that
$Y=G/H$ is a quasi-affine open subvariety of
an affine variety $Y'$.

Then there exists a biholomorphic self-map $\phi$  of $G$
such that the projection map $\pi:G\to G/H$ restricts to an
injective map from $D$ to its image $\phi(D)$ such that $\phi(D)$
is discrete in $Y'$.
\end{proposition}

\begin{proof}
Due to proposition~\ref{torus-quot} we may assume that there is
a maximal torus $T$ of $G$ such that the projection $\rho:G\to G/T$
restricts to an
injective proper map on $D$.
Let $N$ denote the smallest normal subgroup of $G$ containing $T$.
Standard structure theory of linear algebraic groups implies
that $N$ contains one (and therefore all) maximal semisimple
algebraic subgroups of $G$.

We first deal with the case where $H$ does not contain a
maximal semisimple algebraic subgroup of $G$. Then evidently
$N$ is not contained in $H$.
By the definition of $N$ it follows that there exists a conjugate
$gTg^{-1}$ of $T$ which is not contained in $H$. Thus the $T$-action
on $G/H$ is non-trivial and we may use proposition~\ref{pingpong}
to conclude the proof.

Now assume that $H$ contains a maximal semisimple algebraic subgroup.
Let $U$ denote a maximal unipotent algebraic subgroup of $G$.

Clearly $U$ does not contain any semisimple algebraic subgroup.
Therefore, by the preceding consideration, after applying a
biholomorphic self-map of $G$ we may assume that 
the projection $p:G\to G/U$
restricts to an
injective proper map on $D$.
From the Levi decomposition of $G$ we deduce that $G=H\cdot U$.
This implies that the $U$-action on $G/H$ is non-trivial.
Now a second invocation of proposition~\ref{pingpong} yields
the desired result.

Finally we discuss the case where $G$ is unipotent. In this case
$G\simeq\C^n$ as a complex manifold and the projection $G\to G/H$
is simply the projection onto the first $k$ coordinates (for some $0<k<n$)
after a suitable change of coordinates.
Thus the statement follows from the fact that $\phi(D)=\N\times\{0\}^{n-1}$
for some biholomorphic self-map $\phi$.
\end{proof}

\begin{lemma}\label{exists-pingpong}
Let $G$ be a linear-algebraic group with $\dim(G)>1$ which is not
locally isomorphic to $SL_2(\C)$.

Then there exist algebraic subgroups $H$ and $I$ of $G$ such that
\begin{enumerate}
\item
$0<\dim(H),\dim(I)<\dim(G)$.
\item
The quotient varieties $G/H$ and $G/I$ are quasi-affine.
\item
The $H$ action on $G/I$ (by left multiplication) has a dense
orbit.
\end{enumerate}
\end{lemma}
\begin{proof}
Case 1. $G$ is neither reductive nor unipotent. In this case
let $V$ denote the unipotent radical and let $H$ denote a maximal
reductive subgroup. Then $G$ is a semi-direct product of $H$ and $V$.
Therefore $G/V\simeq H$ and $G/H\simeq V$. In particular, both quotient
varieties are affine and $H$ acts transitively on $G/V$..

Case 2. $G$ is unipotent.
Let $G'$ denote the commutator group.
We choose $I$ as a codimension one subgroup of $G$ with $G'\subset I$.
and $H$ as a one-dimensional unipotent subgroup with $H\not\subset I$.
The quotients $G/H$ and $G/I$ are quasi-affine, because $G$ (and therefore
every algebraic subgroup of $G$) is unipotent.

Case 3. $G$ is reductive, but not simple.
Then we can find normal reductive subgroups $H$, $I$ of $G$ such that
$G$ equals the quotient of $I\times H$ by a finite subgroup.
Ther quotients $G/H$ and $G/I$ are affine, because $H$ and $I$
are reductive.

Case 4. $G$ is simple.
We use root system theory. Let $T$ be a ``maximal torus'' of $G$, $\Delta
\subset\Hom(\Lie T,\C^*)$
the root system, $\Delta^+$ a set of positive roots and $S$ a set
of simple positive roots.
Due to our assumption that $G$ is not locally isomorphic to $SL_2(\C)$
we know that the rank of $G$ is at least two. Hence there are at
least two different simple positive roots.
For a simple root $x\in S$ let $\Delta_x^+$ denote the set of
all positive roots which can be written as $\Z^+$-linear combination
of simple roots in $S\setminus\{x\}$.
Furthermore let $T_x=\ker(x:T\to\C^*)$. We choose 
two distinct simple roots $\alpha$, $\beta$.
Associated to $\alpha$, $\beta$ we define two different
algebraic subgroups $H, I$
as follows:
\[
\Lie(H)=\Lie(T_\alpha)+\oplus_{v\in\Delta^+}U_v+\oplus_{v\in\Delta_\alpha^+}
U_{-v}
\]
and 
\[
\Lie(I)=\Lie(T_\beta)+\oplus_{v\in\Delta_\beta^+}U_v+\oplus_{v\in\Delta^+}
U_{-v}
\]
The construction immediately gives us the Levi decompositions of $H$ and
$I$: The Lie algebras given by
\[
\Lie(T_\alpha)+\oplus_{v\in\Delta_\alpha^+}U_v+\oplus_{v\in\Delta_\alpha^+}
U_{-v}
\]
and
\[
\Lie(T_\beta)+\oplus_{v\in\Delta_\beta^+}U_v+\oplus_{v\in\Delta_\beta^+}
U_{-v}
\]
define the Lie algebras of the semisimple subgroups of $G$ corresponding
to the subsets of simple positive roots $S\setminus\{\alpha\}$
and $S\setminus\{\beta\}$ whereas
\[
\oplus_{v \in\Delta^+\setminus\Delta^+_\alpha}U_v
\]
and

\[
\oplus_{v \in\Delta^+\setminus\Delta^+_\beta}U_v
\]
clearly are ideals in the Lie algebra of $H$ resp.~$I$
corresponding to unipotent
subgroups. Thus we see that both $H$ and $I$ are semi-direct products
of a semisimple Lie group with unipotent group.
In particular, neither $H$ nor $I$ admit any
character
(i.e., a homomorphism of algebraic groups to $\C^*$).
It follows
that the quotient varieties $G/H$ and $G/I$ are quasi affine.

Also by construction $\Lie H+\Lie I=\Lie G$.
Hence the $H$-action on $G/I$
has a dense orbit.
Thus $H$ and $I$ have the desired properties.
\end{proof}

\begin{proposition}\label{tame-eq}
Let $G$ be a character-free complex linear algebraic group.

Then every tame discrete set is strongly tame and any two
tame discrete sets are equivalent.
\end{proposition}

\begin{proof}
  It suffices to show: {\em Given two tame discrete sets $D$, $D'$ every
    injective map $\zeta$ from $D$ to $D'$ extends to a biholomorphic
    self-map $F$ of $G$.}

  Upon replacing $D'$ by $\zeta(D)$ it suffices to consider the case
  where $\zeta$ is bijective.
  
  If $G\simeq SL_2(\C)$ or $G\simeq PSL_2(\C)$, this
  is known due to \cite{JW-TAME1}.
  Thus we may assume that $G$ is not locally isomorphic to
  $SL_2(\C)$.

  Fix an exhaustion function $\rho$ on $G$.
  Let $R_k$ be threshold sequence (in the sense of
  Definition~\ref{def-threshold}) for $G$. (The existence is
  due to Theorem~\ref{threshold}.)
  Since $D$ and $D'$ are tame, there are holomorphic self-maps $\phi$
  and $\tilde\phi$ of $G$ such that
  \[
  \#\{ x\in D: \rho(\phi(x))\le R_{2k+1} \} \le k\ \forall k\in\N_0
  \]
  and
  \[
  \#\{ x\in D': \rho(\tilde\phi(x)))\le R_{2k+1} \} \le k\ \forall k\in\N_0.
  \]
  Then
  \[
  \#\{x\in \phi(D)\cup\tilde\phi(D'):
  \rho(x) \le R_{2k+1}\} \le 2k  \ \forall k\in\N_0.
  \]
  It follows that
  \[
  \#\{x\in \phi(D)\cup\tilde\phi(D'):
  \rho(x) \le R_m\} \le m \ \forall m\in\N_0.
  \]
  Since the $(R_k)$ form a threshold sequence, $\phi(D)\cup\tilde\phi(D')$
  is tame. Thus, upon replacing $D$ by $\phi(D)$, $D'$ by $\tilde\phi(D')$
  and $\zeta$ by $\tilde\phi\circ\zeta\phi^{-1}$ we may assume that
  $D\cup D'$ is tame.

  Now we choose $H$ and $I$ as given by Lemma~\ref{exists-pingpong}.
  Then $G/H$ is quasi-affine and there exists an embedding
  $G/H\hookrightarrow V$ into an affine variety $V$.
  We may assume 
  that after application of a suitable biholomorphic
  self-map $\alpha$ of $G$  the natural
  projection $\pi:G \to G/H$ maps
  $D\cup D'$ injectively onto a discrete subset of $V$.
(Proposition~\ref{qa-quot}).
  
  Let $\tau:G\to G/I$ be the natural projection map.
  Proposition~\ref{pingpong}, $(ii)$ implies that
  there exists a biholomorphic self-map $\phi$ of $G$ such that
  \[
  \tau\circ\phi(x)=\tau\circ\phi\circ\zeta(x)\ \forall x\in D.
  \]
  and such that $\tau\circ\phi$ is injective on $D'$.

  It follows that there is a map $\beta:D\to I$
  such that
  \[
  \phi(x)\beta(x)=\phi(\zeta(x))\ \forall x\in D
  \]

  In addition, due to Proposition~\ref{pingpong}
  $\tau(\phi(D'))$ is discrete in some affine completion $Y'$ of $Y=G/I$.
  Since $\zeta$ is bijective, $\tau\phi$ also maps $D$ bijectively
  onto a discrete subset of the affine variety $Y'$.
  
  Hence we may choose a holomorphic map $\gamma:G/I\to H$ with
  \[
  \gamma(\tau(\phi(x)))=\beta(x)\ \forall x\in D
  \]
  and define a self-map $\Psi$ of $G$ as
  \[
  \Psi:y\mapsto y\cdot \gamma(\tau(y))
  \]
  Then
  $y\mapsto \phi^{-1}\circ\Psi\circ\phi$ is a biholomorphic
  self-map of $G$ whose restriction to $D$ agrees with $\zeta$.
  \end{proof}

\begin{proposition}\label{sub-uni}
Let $G$ be a complex linear algebraic group
and let $H\subsetneq G$ be a unipotent algebraic subgroup.

Let $D\subset H$ be a discrete subset. 

Then $D$ is a tame discrete subset of $G$.
\end{proposition}

\begin{proof}
First, we assume that $G$ is not unipotent. Let $T$ be a ``maximal
torus''. Note that $\dim(T)>0$, because $G$ is not unipotent.
For any $g\in G$ the group $gHg^{-1}$ is unipotent and therefore
has trivial intersection with $T$. It follows that the $H$-action
on $G/T$ is free. In particular, all orbits have the same dimension.
Since the group action is algebraic, it follows that all orbits are
closed. Thus $H$ acts freely on $G/T$ with closed orbits.
Furthermore $G/T$ is an affine variety, because $T$ is reductive.
Thus we may invoke Proposition~\ref{proj-tame} to conclude that $D$ is tame.

Second, assume that $G$ is unipotent.
Now $H$ is a non-trivial unipotent subgroup of a unipotent group.
By standard theory of nilpotent Lie groups it follows
that $H\cdot G'\ne G$ where $G'$ denotes the commutator group of $G$.
Therefore there exists a complex one-parameter Lie subgroup $A$ of $G$
with $A\not\subset HG'$.
Because $G$ is unipotent, its exponential map is algebraic and every
connected complex Lie subgroup of $G$ is algebraic. Thus $A$ is an
algebraic subgroup.
Now $gHg^{-1}$ is contained in $HG'$ for every $g\in G$, because
$HG'$ is normal in $G$. It follows that $gHg^{-1}\cap A$ is trivial for
every $g\in G$, i.e., $H$ acts freely on $G/A$. By the same reasoning
as in the first case above we know that every $H$-orbit in $G/T$ is closed.
Finally note that
$G/A$ is isomorphic to an affine space, because both $G$ and $A$ are
unipotent.
\end{proof}

\section{Unions of tame sets}

\begin{lemma}
Let $G$ be a complex linear algebraic group with maximal 
connected reductive commutative subgroup (``maximal torus'') $T$.

If $\dim(T)>1$, then $T$ contains two subgroups isomorphic to $\C^*$
which are not conjugate in $G$.
\end{lemma}

\begin{proof}
Let $R$ denote the radical of $G$. 
If $T\cap R=T$, then $T=G$, i.e. $G$ is solvable.
Then $T$ maps injectively into $G/[G:G]$, which implies that two
different subgroups of $T$ can not be conjugate.

If $0<\dim(R\cap T)<\dim T$, the assertion is obvious, because
subgroups of $R\cap T$
can  not be conjugate to subgroups which are not contained in $R$.

The remaining case is the case where $\dim(R\cap T)=0$. Then $R$
is unipotent, and there is no loss in generality in replacing $G$ by
$G/R$. Hence we may assume that $G$ is semisimple.
Let $\Delta \subset \Hom(T,\C^*)$ be the root system, i.e.,
\[
\Delta=\left\{  \lambda\in  \Hom(T,\C^*):\exists v\in \Lie(G)\setminus\{0\}:
ad(g)(v)=(\lambda(g))v\ \forall g\in T\right\}
\]
We define
\[
R_\alpha=\{v\in \Lie(G):
ad(g)(v)=(\lambda(g))v\ \forall g\in T\}
\]
for $\alpha\in \Delta$.

Now $\Sigma=\cup_{\alpha\in\Delta} \ker\alpha$ denotes the set of all
$g\in T$ whose centralizer contains a unipotent element.
Evidently a one-dimensional subgroup contained in $\Sigma$ can not
be conjugate to  a subgroup of $T$ which is not contained in $\Sigma$.
Therefore there do exist non-conjugate subgroups isomorphic to $\C^*$.
\end{proof}

\begin{proposition}\label{tech-union}
Let $G$ be a connected complex Lie group with 
positive-dimensional  connected Lie subgroups 
$H$ and $I$ such that $G/H$ and $G/I$ are Stein manifolds
 and such that
the diagonal map $\delta:G\to (G/H)\times (G/I)$ is a 
proper map.

Then every discrete subset of $G$ is the union of two tame
subsets.
\end{proposition}

\begin{proof}
Let $\rho_1$ resp.~$\rho_2$ be an exhaustion function on
$G/H$ and $G/I$. Via $\delta$ we obtain an induced
exhaustion function $\rho$ on $G$:
\[
\rho(g)=\rho_1(gH)+\rho_2(gI)
\]
Let $D$ be a discrete subset of $G$. We define
$D_1=\{g\in D:\rho_1(gH)\ge \rho_2(gI)\}$
and
$D_2=\{g\in D:\rho_1(gH) < \rho_2(gI)\}$
We observe that $\rho(g)\le 2\rho_1(gH)$ for every $g\in D_1$.
Hence $D_1$ admits a proper map to $G/H$. This implies
that $D_1$ is tame due to Proposition~\ref{proj-tame}.

Similarily for $D_2$.
\end{proof}

\begin{theorem}\label{union-2-tame}
  Let $G$ be a connected complex linear algebraic group,
  $\dim(G)\ge 2$.

Then every discrete subset of $G$ is the union of two  tame
subsets.
\end{theorem}

\begin{proof}
Let $T$ be a maximal torus.

Case 1: $\dim(T)>1$. Then $G$ admits two subgroups $H$, $I$ which are
isomorphic to $\C^*$ and which are not conjugate.
The quotients $G/H$ and $G/I$ are Stein, because $H$ and $I$
are reductive.
We consider the $G$-action on $Y=(G/H)\times(G/I)$. The isotropy group
at $(xH,yI)$ equals $xHx^{-1}\cap yIy^{-1}$. This is necessarily finite,
since $H$ and $I$ are one-dimensional and not conjugate.
Thus all isotropy groups are finite. It follows that each $G$-orbit in $Y$
has the same dimension. By algebraicity this implies that every
$G$-orbit in $Y$ is closed. Thus every orbit map $G\to Y$
(with $g\mapsto (gxH,gyI)$) is a proper map.
Hence we may apply
Proposition~\ref{tech-union}
and deduce the desired assertion.

Case 2: $\dim(T)<2$ and $G$ is solvable. If $G$ is solvable,
$G$ is biholomorphic to $(\C)^k\times(\C^*)^m$.
Thus $G$ is isomorphic as a variety to a product $\tilde G$ of two linear
algebraic groups $H$ and $I$, which allows us to employ
Proposition~\ref{tech-union}.

Case 3:
$\dim(T)<2$ and $G$ is neither solvable nor semisimple.
Then necessarily $T$ is one-dimensional and a maximal connected
semisimple subgroup $S$ of $G$ is isomorphic to $SL_2(\C)$.
Here we may argue with $H=R$ and $I=T$.

Case 4:
$G$ is (up to a finite covering) isomorphic to $SL_2(\C)$.
Here we may refer to our earlier paper \cite{JW-TAME1}.
\end{proof}

\section{Oka theory}

A complex manifold $X$ is called an {\em Oka manifold} if for
every compact convex subset $C$ of $\C^n$ the space of
holomorphic functions $Hol(\C^n,X)$ is dense
(with respect to the topology of locally uniform convergence)
in
$Hol(C,X)$. This is known to imply the property
every Stein complex manifold $Y$ and every continuous map from
$Y$ to $X$ there exists a homotopic holomorphic map.

One instrument to prove this Oka property is a ``spray'':
\begin{definition}
  Let $X$ be a complex manifold. A spray is a pair $(E,\pi,s)$ where
  $\pi:E\to X$ is a  vector bundle and $s:E\to X$ is holomorphic map
  such that
    $s=\pi$ on the zero-section $S$ of $E$.

  A spray is called ``dominating''if
    $s|_{E_x}:E_x\to X$ has maximal rank at $x$ for every $x\in S$.
\end{definition}

A complex manifold with a dominating
spray in this sense is necessarily an Oka
manifold.

This generalizes to the notion of ``subellipticity'' or the notion
of a dominating family of sprays.

\begin{theorem}\label{subelliptic}
  Let $X$ be a complex manifold and let $\pi:E_i\to X$ be a finite family
  of holomorphic vector bundles with maps $s_i:E_i\to X$ such that
  $s_i$ agrees with $\pi_i$ on the zero-section of $E_i$ and such
  that for each point $x\in X$ we have $T_xX=\oplus_k (ds_i)(T_{p_i}E_i)$
  where $p_i$ is the point in the zero section of $E_i$ above $x$.

  Then $X$ is a Oka manifold.
\end{theorem}
  
See \cite{Forstneric} for more information about Oka manifolds.

\begin{proposition}\label{oka-tool}
Let $X$ be a complex manifold on which its automorphism group
$\Aut(X)$ acts transitively. Let $\pi:E\to X$ be a vector bundle,
and let $\tau:E\to X$ be a holomorphic map such that
$\tau$ equals the identity map  on the zero-section and such that
there exists a point where $\tau$ has maximal rank.

Then $X$ is an Oka manifold.
\end{proposition}

\begin{proof}
Let $S$ denote the zero-section of $E$. Let $Z_0=Z$ denote the set of
points of $S\simeq X$ where $\tau$ does not has maximal rank.
For each $p\in Z$ the set $\{\phi\in Aut(X):\phi(p)\in Z_0\}$ is nowhere
dense. We choose a countable subset $D\subset Z_0$ such that $D$
intersects every irreducible component of $Z_0$.
Due to Corollary~\ref{prep-baire},
there is an element $g_0\in Aut(X)$ with
$g_0(D)\subset X\setminus Z$.
Define $Z_1=Z_0\cap g_0(Z_0)$. Then $\dim(Z_1)<\dim(Z_0)$.
Continuing recursively, we obtain finitely many elements $g_0,\ldots,g_s$
of $Aut(X)$ such that $\cap_k g_k(Z_k)$ is empty.

Now the pull-back bundles $g_k^*E$ combined with the maps
$v\mapsto g_k(\tau(v))$ constitute a subelliptic family of
dominating sprays in the sense of Theorem~\ref{subelliptic}.

Therefore $X$ is an Oka manifold.
\end{proof}

\begin{lemma}\label{8.2}
Let $V$ be a normal quasi affine variety, and let $f:\C\to V$ be an injective
morphism. Then $f$ is a proper map, and every holomorphic function
on $f(\C)$ extends to a holomorphic function on $V$.
\end{lemma}

\begin{proof}
Let $V$ be embedded into a normal affine variety $V_1$ which in turn is embedded
into a normal projective variety $\bar V$.
Then $f$ extends to a morphism from $\P_1=\C\cup\{\infty\}$ to $\bar V$.
Since $f(\P_1)$ is compact and $V_1$ is affine, it is clear that
$f(\P_1)\not\subset V_1$. Therefore $f(\infty)\not\in V_1$. It follows
that $f(\C)=V_1\cap f(\P_1)$ which in turn implies that $f(\C)$ is closed
in $V_1$. Combined with $V_1$ being normal and the injectivity of $f$,
we deduce that $f:\C\to V_1$ is a closed embedding and in particular 
proper.

Now $f(\C)$ is a closed analytic subset of the Stein space $V_1$. Therefore
every holomorphic function extends from $f(\C)$ to $V_1$.
Now the statement of the lemma is obvious, since $V$ is an open subset
of $V_1$. 
\end{proof}

\begin{proposition}
Let $G$ be a complex linear algebraic group, $U$ a one-di\-men\-si\-onal unipotent
subgroup, 
and $H$ a connected algebraic
subgroup of $G$ with $H\cap U=\{e\}$.

If $H$ is either reductive or unipotent, then every holomorphic
function on $U$ extends to an $H$-invariant holomorphic function
on $G$. 
\end{proposition}

\begin{proof}
  Since $H$ is reductive or unipotent, we know that $G/H$ is quasi-affine.
  Due to $H\cap U=\{e\}$ the projection map $\pi:G\to G/H$ restricts to
  an injective map on $U$. Lemma \ref{8.2} implies that the image
  $\pi(U)$ is a closed subvariety of $G/H$ such that every holomorphic
  function on $\pi(U)$ extends to the complex manifold $G/H$.
  Since holomorphic functions on $G/H$ correspond to $H$-invariant
  holomorphic functions on $G$, the assertion follows.
\end{proof}

\begin{theorem}\label{oka1}
  Let $G$ be a complex linear algebraic group, and let $D$ be
  a tame discrete subset.

Then $X=G\setminus D$ is a Oka manifold.

Furthermore there exists an infinite discrete subset $D'\subset G$ such
that $G\setminus D'$ is not an Oka manifold.
\end{theorem}

\begin{proof}
  Let $U$ be a one-dimensional unipotent subgroup of $G$, i.e., an algebraic
  subgroup isomorphic to the additive group $(\C,+)$. (Such a subgroup
  exists due to Lemma~\ref{exist-unip}.)
  Thanks to the equivalence $(i)\iff(iv)$ of Theorem~\ref{main-1} we
  may assume that $D\subset U$.

  Let $N$ denote the ``nil cone'' of $G$  and 
let $v=(v_1,\ldots,v_n)$ be a vector space basis of $\Lie(G)$ such that
$v_i\in N$ for all $i$. (This is possible by Lemma~\ref{nil-cone}.)
Since $ad(v_i)$ is nilpotent,
the algebraic groups $U_i$ generated by $\{\exp (v_it):t\in\C\}$ are
all one-dimensional. We choose $v$ such that every $U_i$ is
different from $U$.

We define
\[
\Sigma_v=\cup_{1\le i\le n}\cup_{x\in D}U_i\cdot x.
\]
This is a closed analytic subset of $G$,
because $D$ is discrete and the orbits
of the $U_i$ are closed in $G$.

We obtain flows $\phi_i$:
\[
\phi_i(t;g)\stackrel{def}=\exp(v_it)\cdot g
\]

We choose $U_i$-invariant holomorphic functions $\alpha_i$ on $G$
such that $\{g\in U:\alpha_i(g)=0\}=D$ for every $i$.
(This is possible due to Lemma~\ref{8.2}.)

Now we can modify the flows $\phi_i$:
\[
\tilde\phi_i(t;g)= \exp(v_it\alpha_i(g))\cdot g
\]

By construction, $\tilde\phi_i$ acts trivially on
\[
\{u\cdot x: u\in U_i, x\in D\}
\]

We define further a map $\Psi_v:\C^n\times G\to G$ as follows
\[
\Psi_v(t_1,\ldots,t_n;p)=
\tilde\phi_1(t_1)\circ\ldots\circ\tilde\phi_n(t_n)(p)
\]
Since the $v_1,\ldots,v_n$ form a basis of $\Lie(G)$,
it is clear that for a fixed element $p\in G$ the map
$t\mapsto \Psi_v(t)$ has maximal rank at $(0;p)$ if 
$\alpha_i(p)\ne 0\ \forall i$.

Let $A$ denote the set of all holomorphic automorphisms of $G$
which fix $D$ point-wise.
We recall that $A$ acts transitively on $Y=G\setminus D$.

Therefore we may now invoke Proposition~\ref{oka-tool}
(with the trivial vector bundle as $E$) and deduce that
$G\setminus D$ is an Oka manifold.

Finally we recall that every Stein manifold (hence in particular
every complex linear algebraic group) $G$ admits
a discrete subset $D'$
such that every holomorphic map from $\C^n$ ($n=\dim_C(X)$)
to $G\setminus D'$ is
degenerate
(\cite{JW-LDS}).
Evidently it follows that the
complement $G\setminus D'$ is not an Oka manifold.
\end{proof}

\section{Oka of blow up}

\begin{lemma}\label{oka-cn}
  Let $D=\Z\times\{0\}^{n-1}\subset\C^n$.

  Then the manifold $\hat X$ obtained from $\C^n$ by blowing up $D$ is Oka.
\end{lemma}

\begin{proof}
  We observe that
  \[
  D=\{(z,0\ldots,0):f(z)=0\}
  \]
  for $f(z)=\exp(2\pi i z)-1$.
  The blow up of $\C^n$ along $D$ may be described as
  \begin{multline*}
  \hat X=\{(z_1,\ldots,z_n;[w_1:\ldots:w_n])\in \C^n\times\P_{n-1}:\\
  f(z_1)w_j=z_jw_1 \ \forall j>1, z_jw_k=z_kw_j\ \forall j,k>1\}
  \end{multline*}
  Now 
  \[
  \phi:(z_1,\ldots,z_n]\mapsto
  \left(z_1,f(z_1)z_2,\ldots,f(z_1)z_n;[1:z_2:\ldots:z_n]\right)
  \]
  defines a biholomorphic map from $\C^n$ to $\hat X\setminus H$
  where 
  \[
  H=\{
  (z_1,\ldots,z_n;[w_1:\ldots:w_n]):f(z_1)=0=w_1\}
  \]
  is the strict transform of the hypersurface defined by $f(z_1)=0$.
  For every linear map $L:\C^{n-1}\to\C$ we obtain an automorphism
  $\alpha_L$ of
  $\C^n$ given by
  \[
  (z_1,z_2,\ldots,z_n)\mapsto (z_1+L(z_2,\ldots,z_n),z_2,\ldots,z_n).
  \]
  This automorphism of $\C^n$ lifts to an automorphism $\beta_L$ of
  $\hat X$,
  because it stabilizes the exceptional locus of the blow up defined
  by $f(z_1)=0=z_2=\ldots=z_n$.
  Now $\beta_L\circ\phi$ defines a biholomorphic from $\C^n$ to the
  complement of the strict transform of the hypersurface defined
  by $f(z_1+L(z_2,\ldots,z_n))=0$.
  By taking $(n-1)$ linear independent linear forms $L_j$ and
  the corresponding maps $\beta_{L_j}\circ \phi$ we see that $\hat X$
  may be covered by finitely many Zariski open subsets biholomorphic to
  $\C^n$. Since $\C^n$ is Oka, the Oka property for $\hat X$ now follows
  from \cite{YK}.
\end{proof}
\begin{lemma}\label{blow-up-solv}
Let $H=(\C^*)^k\times\C^n$ with $k\ge 0$, $n\ge 1$
and $k+n\ge 2$. Let $D$ be a discrete subset of
$H$ which is contained in $(\C^*)^k\times\{\C\}^{n-1}\times\{0\}$
and  let $Z$ denote the complex manifold obtained
from $H$ by blowing up $D$. 

Then $Z$ is a Oka manifold.
\end{lemma}

\begin{proof}
Let $\tau:\tilde H=\C^{k+n}\to H$ denote the universal covering
of $H$. Then the universal
covering of $Z$ is the complex manifold obtained from $\tilde H$
by blowing up the discrete subset $D'=\tau^{-1}(D)$.
Observe that $D'$ is tame in $\tilde H$, 
because $D'\subset\C^{k+n-1}\times\{0\}$.
(\cite{RR}, Corollary 3.6.a.)
Hence $\tilde Z$ is biholomorphic to the complex manifold
obtained from $\C^{k+n}$ by blowing up $\Z\times\{0\}^{n+k-1}$.
Therefore $\tilde Z$ is Oka. 
(Lemma \ref{oka-cn}).
Since $\tilde Z\to Z$ is a holomorphic covering,
$Z$ must be Oka, too (\cite{Forstneric}, Proposition 5.5.2).
\end{proof}

\begin{proposition}\label{oka2}
Let $G$ a character-free complex linear algebraic group.
Let $D$ be a tame discrete subset and let $X$ denote the complex manifold
obtained from $G$ by blowing up $D$.

Then
$X$ is a Oka manifold.
\end{proposition}

\begin{proof}
Let $S$ denote a maximal semisimple Lie subgroup of $G$ and let $R$
denote the unipotent radical. Let $T$ be a maximal torus of $S$.

 Case $(1)$ Assume that $S=\{e\}$. Then $G=R$. Hence
$G$ is biholomorphic to some $\C^n$ with $n\ge 2$  and the assertion 
follows from lemma~\ref{blow-up-solv}.

Case $(2)$: Assume that $\dim(S)>0$. Then $T$ is non-trivial and thus
unbounded. Consequently
$T$ contains a tame discrete subset $D'$ (Corollary~\ref{unbounded}).
Since all tame discrete subsets of $G$ are equivalent,
we may assume that $D\subset T$ (Corollary~\ref{cor-equiv}).

Fix a root system for $S$ and let $U^+$ resp.~$U^-$ denote the maximal
unipotent subgroup of $D$ corresponding to all positive resp.~
negative roots. Next we claim that the map
$\alpha:W_0=U^-\times T\times U^+\to S$
defined as
\[
(x,y,z)\mapsto x\cdot y\cdot z
\]
is an open embedding.
It suffice to check injectivity. However
\[
x\cdot y\cdot z = x'\cdot y'\cdot z',\ x,x'\in U^-, y,y'\in T, z,z'\in U^+
\]
implies
\[
(x')^{-1}\cdot x =  y'\cdot z'\cdot  (y\cdot z')^{-1}
\]
Now $(x')^{-1}\cdot x $ is an element of the unipotent group $U^-$
while the right hand side is an element in the Borel group $TU^+$.
It follows that both sides equal $e_G$ which implies the claim.

Now $G$ is a semi-direct product of $S$ with the unipotent radical $R$,
thus we obtain a Zariski open embedding of $W=W_0\times R$ into
$G$ via $(w,r)\mapsto \alpha(w)\cdot r$.

The complex manifold $\hat W$ obtained from $W$ by blowing up $D$
(remember that we assumed $D\subset T\subset W$ is Oka due to
lemma~\ref{blow-up-solv}.

Define $V=(S\setminus T)\cdot R$.
Now $G$ is flexible
(see \cite{FKZ})
and $T$ is of codimension at least $2$ in $G$. Therefore
$V$ is still flexible
(\cite{FKZ}, Theorem 1.1) and as a consequence Oka
(\cite{Forstneric}, Proposition 5.5.17).

We thus have found two Zariski open submanifolds of $X$, namely $V$
and $\hat W$ which are both Oka.

Therefore $X$ is Oka
due to \cite{YK}, Theorem 1.4.
\end{proof}

\begin{remark}
  (\cite{JW-LDS}). Evidently $X\setminus D$ is not Oka.
\end{remark}

\section{Generalized Overshears}

\begin{proposition}\label{prescribe-auto}
  Let $Y_1$ be a Stein space, $E$ a nowhere dense
  analytic subset of $Y_1$,
$Y=Y_1\setminus E$, $\pi:X\to Y$ a $\C$-principal bundle with a meromorphic
section $\sigma$, and $D\subset Y_1$ a discrete subset such that
$D\cap(E\cup\{\sigma=\infty\})=\{\}$.
For each $p\in D$, let an automorphism $\alpha_p$ of the fiber
$X_p=\pi^{-1}(p)$ be given.

Then there exists an automorphism $\phi$ of $X$ such that $\pi=\pi\circ\phi$
and such that $\phi|_{X_p}=\alpha_p$ for all $p\in D$.
\end{proposition}
\begin{proof}
  Let $(U_i)_i$ be a trivializing open cover of $Y$ and
  let $\phi_{ij}:U_i\cap U_j\to\C$ be
  the corresponding transition functions. Then
  the meromorphic section $\sigma$ is given by a
collection of meromorphic functions $s_i\in \M(U_i)$ such that
$s_i=s_j+\phi_{ij}$ on each $U_i\cap U_j$.
A fiber-preserving automorphism $\phi$ of $X$ is given
by a holomorphic function $\lambda:Y_1\to\C^*$ and a family of
holomorphic functions $c_i\in\O(U_i)$ such that
in local coordinates $\phi$  may be described as
\[
\phi: (x,z_i)\mapsto (x,\lambda(x) z_i+c_i(x).
\]
The functions $c_i$ must fulfill
\[
\lambda z_j+c_j+\phi_{ij} = \lambda (z_j+\phi_{ij})+c_i,
\]
i.e., $c_i=c_j+\phi_{ij}(1-\lambda)$.

Now we may choose $\lambda\in\O^*(Y_1)$ arbitrarily subject only to
the condition that $1-\lambda$ has to vanish
(with the respective multiplicities)
along the poles of $\sigma$.
Once $\lambda$ is fixed, we may choose $h\in\O(Y_1)$ freely,
and define
\[
c_i=s_i(1-\lambda)+h.
\]
This yields the assertion, because on a Stein space we may prescribe the
values of holomorphic functions on a discrete subset.
\end{proof}

\begin{proposition}
Let $G$ be a complex linear algebraic group, $H$ 
a commutative reductive subgroup,
$U\simeq (\C,+)$ a one-dimensional unipotent subgroup and let
$D\subset G$ be a discrete subset.

Let $j:G/U=Y\hookrightarrow Y_1$ be an open algebraic embedding into
an affine variety $Y_1$, let $\pi:G\to G/U$, $\tau:G\to G/H$ 
denote the natural
projections and assume that $\pi(D)$ is discrete in $Y_1$.

Then there exists a holomorphic automorphism $\phi$ of $G$
such that $\tau\circ\phi|_D:D\to G/H$ is a proper map.
\end{proposition}

\begin{proof}
There is no non-constant algebraic morphism from $\C$ to $H$, therefore
the restriction of $\tau$ to an $U$-orbit $gU$ is always non-constant.
Since $G/H$ is an affine variety
(because $H$ is reductive),
it follows that the morphism
$\tau|_{gU}:gU\to G/H$ is actually a finite proper morphism for
every $g$.

We fix a proper exhaustion function $\rho$ on $G/H$ and
a bijective map $\zeta:\N\to \pi(D)$.

Since $G\to G/U$ is an algebraic $\C$-bundle, there exists a meromorphic
section $\sigma$. Let $P\subset Y=G/U$ denote its polar divisor.
Since $D$ is countable, we may replace $D$ by a translate and thereby
without loss of generality assume that $\pi(D)$ does not intersect $P$.

Therefore we may, using the above Proposition~\ref{prescribe-auto},
prescribe automorphisms
along $\pi(D)$. 

From this fact we deduce the assertion in the following way:

We fix an enumeration $\pi(D)=\{q_k:k\in\N\}$.
We choose an exhaustion function $\rho $ on $G/H$.
We observe: {\em given a discrete subset $\Lambda$ of $\C$ and a
  compact subset $C$ of $\C$, its is always possible to choose an
  automorphism $\zeta:z\mapsto\alpha z+\beta$ ($\alpha\in\C^*$, $\beta\in\C$)
  such that $\zeta(\Lambda)\cap K\ne\{\}$.)}

Therefore, for each $n\in\N$, we may choose an automorphism $\zeta_k$
of $L_k\stackrel{def}{=}\pi^{-1}(q_k)$ such that
\[
\tau(\zeta_k(D\cap L_k))\ \cap \ \{x\in G/H; \rho(x)\le n\}=\{\}.
\]
Now Proposition~\ref{prescribe-auto} yields an automorphism $\phi$
of the complex manifold $G$ which for each $k\in\N$ agrees with $\zeta_k$
on $L_k$.
This implies that
\[
\{x\in D: \rho(\tau(\phi(x)))\le n\}
\]
is finite for every $n\in\N$.
\end{proof}

\begin{corollary}\label{overshear-tame}
Let $G$ be a complex linear algebraic group,
$U\simeq (\C,+)$ a one-dimensional unipotent subgroup and let
$D\subset G$ be a discrete subset.

Let $j:G/U=Y\hookrightarrow Y_1$ be an open algebraic embedding into
an affine variety $Y_1$, let $\pi:G\to G/U$, denote the natural
projection and assume that $\pi(D)$ is discrete in $Y_1$.

Then $D$ is tame.
\end{corollary}

\begin{proof}
  If $G$ is unipotent, then $G\simeq\C^n$ as a complex manifold. Moreover
  in this case we may choose coordinates on $\C^n$ such that $\pi$ is a
  linear projection. Then the assertion follows from standard facts on
  tame discrete sets in $\C^n$. (\cite{RR}.)

  Hence we may assume that $G$ is not unipotent. Then there exists a
  positive-dimensional commutative reductive subgroup of $G$ and we may
  invoke the preceding proposition. In Combination with
  Proposition~\ref{proj-tame} we obtain tameness of $D$.
\end{proof}

\section{Arithmeticity.}

In this setion we exploit aithmeticity conditions to show that certain
dsicrete subgroups are tame as discrete subsets.

The fact that $\Z$ is discrete in $\R$ generalizes as follows.
Let $K$ be a number field of degree $d=2s+r$ with $r$ real embeddings
$\tau_i:K\to\R$ and $s$ pairs of conjugate complex embeddings $(\sigma_j,
 \bar\sigma_j)$,
$\sigma_j:K\to\C$.
Let $\Ok$ be its ring of algebraic integers, i.e., the integral closure
of $\Z$ in $K$.
Then the diagonal map
$x\mapsto (\tau_1(x),\ldots,\tau_r(x);
\sigma_1(x),\ldots,\sigma_s(x))$ embedds
$\Ok$ as a discrete subring in  $\R^r\times\C^s$.
If $G$ is an algebraic group defined over $K$, in the same way
$G(\Ok)$ embedds as a discrete subgroup in $G(\R)^r\times G(\C)^s$.
 
\begin{fact}
Let $K$ be a number field with ring of algebraic integers $\Ok$.
Let $P:\A^n\to\A^m$ be a polynomial map which is defined over $K$,
i.e., all the coefficients of the polynomial functions $P_1,\ldots, P_m$
are in $K$. 

Then $P(\Ok^n)$ is contained in $\frac{1}{N}\Ok^m$ for some $N\in\N$.
\end{fact}

In fact, we may simply chose $N$ as the lowest common multiple of
the denominators of the coefficients of $P$.

This implies the following fact which is crucial in the sequel.

\begin{fact}\label{fact-arith}
Let $X, Y$ be affine algebraic varieties defined over a number field $K$ and
let $f:X\to Y$ be a $K$-morphism. Fix a closed embedding of $X$ into
an affine space $A^n$.

Then there is a closed embedding  of $Y$ into an affine space $\A^m$
such that $f(X\cap\A^n(\O_k))\subset Y\cap\A^m(\O_k)$.
\end{fact}

Given a group $G$ with two subgroups $\Gamma_i$ $i=1,2$) these two
subgroups are called {\em commensurable} if the intersection
$\Gamma_1\cap\Gamma_2$ is of finite index in both $\gamma_1$ and $\Gamma_2$.

\begin{theorem}\label{arithm}
Let $K$ be a number field of degree $d=2s+r$ with $r$ real embeddings
$\tau_i:K\to\R$ and $s$ pairs of conjugate complex embeddings $(\sigma_j,
 \bar\sigma_j)$,
$\sigma_j:K\to\C$.
 Let $\Ok$ be its ring of algebraic integers.
  Let $G$ be a character-free complex linear algebraic $K$-group
    such that there is a unipotent element $u\in G(K)\setminus\{e\}$.
    Let $H=G(\C)^{r+s}$ and let $G(\Ok)$ be embedded into $H$ via the diagonal
    embedding of $G(K)$ into $G(\R)^r\times G(\C)^s\subset G(\C)^{r+s})$.
    Let $\Gamma$ be a discrete subgroup which is commensurable to $G(\O)$.

  Then $\Gamma$ is a tame discrete subset of the complex manifold $H$.
\end{theorem}
\begin{proof}
  Let $U$ denote the Zariski closure of $\{u^n:n\in\Z\}$. Since $u$ is a $K$-rational
  unipotent element, $U$ is a one-dimensional unipotent $K$-group.
  Consider the
  quotient morphism: $\tau:H\to H/U=V_0$.
  Note that $V_0$ is a quasi affine variety defined
  over $K$. Let $j:V_0\to V$ be a dense open embedding into an affine $K$-variety $V$.
  Observe that $j(\tau(G(\O)))$ is a discrete subset of $V(\C)$
  (cf.~\ref{fact-arith}).
  Since $\Gamma$ and $G(\Ok)$ are commensurable, there exists a finite
  subset $F$ of $H$ with $\Gamma\subset F\cdot G(\Ok)$.
  Therefore $j(\tau(\Gamma)))\subset F\cdot\tau(G(\Ok))$
  is likewise discrete in $V(\C)$.
  Hence tameness of $\Gamma$  follows from Corollary~\ref{overshear-tame}.
\end{proof}

\section{Lattices}
A {\em lattice} $\Gamma$ (in a Lie group $G$) is a discrete subgroup such that
$G/\Gamma$ has finite volume (with respect to the Haar measure).
Lattices are also called discrete subgroups {\em of finite covolume}.
Such a lattice $\Gamma$
is called {\em cocompact} or {\em uniform} if the quotient
$G/\Gamma$
is compact.

A lattice $\Gamma$ in a semisimple complex algebraic group $S$ is called
{\em irreducible} if there is no product decomposition $S=S_1\cdot S_2$
with
\begin{enumerate}
\item
  Both $S_i$ are normal semisimple subgroups of $S$ with
  $0<\dim(S_i)<\dim S$,
\item
  $S_1\cap S_2$ is finite and
\item
  $\Gamma$ and $(\Gamma\cap S_1)\cdot(\Gamma\cap S_2)$ are
  commensurable.
\end{enumerate}

\begin{theorem}
  Let $\Gamma$ be a
  non-cocompact lattice in a semisimple complex
  linear algebraic group $G$.

  Then $\Gamma$ is a tame discrete subset of the complex manifold $G$.
\end{theorem}
\begin{proof}
  The general case is easily reduced to the ace
  where $\Gamma$ is an irreducible lattice.
  Hence from now on  we  assume that $\Gamma$ is irreducible.
  
  Next we recall that due to the theorem of Kazhdan-Margulis $\Gamma$
  contains a unipotent element. (see \cite{KM})
  
  If $G$ is of ``rank one'', i.e., locally isomorphic to $SL_2(\C)$,
  tameness of $\Gamma$ follows from
  Proposition~\ref{SL2-unipotent} below.

  Thus we may assume that $G$ has rank at least two. Now we can apply
  the arithmeticity results of Margulis (see \cite{Margulis}).
  This implies that there is a number field $K$, a $K$-group
  $H$, a choice of complex embeddings $\sigma_i$ of $K$ and
  real embeddings $\tau_j$ such that $\Gamma$ is commensurable to
  $\pi(\delta(H(\Ok))$
  where $\delta$ is the diagonal embedding of $H$ into
  \[S=\left( \Pi_i{}^{\sigma_i}H(\C)\right)\times
  \left( \Pi_i{}^{\sigma_i}H(\R)\right)
  \]
  and $\pi$ denotes the quotient by the compact group
  \[
  U=\left( \Pi_i{}^{\sigma_i}H(\R)\right)
  \]
  However, since $\Gamma$ and therefore also   $\pi(\delta(H(\Ok))$
  is non-cocompact, there is a unipotent element.
  If $H(\R)$ is compact, it can not contain a unipotent element.
  Thus in our case the compact group $U$ must be trivial.
  This allows us to invoke
  theorem~\ref{arithm}
  in order to deduce tameness of $\Gamma$.
\end{proof}

  \section{Discrete subgroups of $SL_2(\C)$}
  
\begin{proposition}\label{SL2-unipotent}
Let $\Gamma $ be a discrete subgroup of $SL_2(\C)$ or $PSL_2(\C)$.

If $\Gamma\setminus\{e\}$ contains a unipotent element, then
$\Gamma$ is a tame discrete subset.
\end{proposition}

\begin{proof}
For simplicity we consider only the case where the group is $SL_2(\C)$.
Let $u$ be a unipotent element of $\Gamma\setminus\{e\}$.
All unipotent elements of $SL_2(\C)$ are conjugate, hence
we may without loss of generality assume that
\[
u=\begin{pmatrix}
1 & 1 \\
0 & 1 \\
\end{pmatrix}.
\]
Let $\pi:SL_2(\C)\to\C^2$ denote the projection onto the first column.
We claim that the image $\pi(\Gamma)$ is discrete in $\C^2$.
To verify this, assume to the contrary that there exists a sequence
$\gamma_n\in\Gamma$ such that $\pi(\gamma_n)$ is convergent in $\C^2$
with $\lim\pi(\gamma_n)=v$, but $\pi(\gamma_n)\ne v$.
Let
\[
\gamma_n=\begin{pmatrix}
a_n & b_n \\
c_n & d_n\\
\end{pmatrix}.
\]
Our assumptions mean that $\lim (a_n,c_n)=(a,c)$ for some
$(a,c)\in\C^2$, but
$(a_n,c_n)\ne (a,c)$ for infinitely many $n\in\N$.

We calculate $\gamma_n\cdot u \cdot \gamma_n^{-1}$:
\[
\gamma_n\cdot u \cdot \gamma_n^{-1}
=
\begin{pmatrix}
a_n & b_n \\
c_n & d_n\\
\end{pmatrix}
\cdot 
\begin{pmatrix}
1 & 1 \\
0 & 1 \\
\end{pmatrix}
\cdot
\begin{pmatrix}
d_n & -b_n \\
-c_n & a_n\\
\end{pmatrix}
=
\begin{pmatrix}
1 -a_nc_n & a_n^2 \\
- c_n^2 & 1+a_nc_n \\
\end{pmatrix}
\]
(Here we use that the determinant $\det(\gamma_n)=a_nd_n-b_nc_n$ equals 
$1$.)
It follows that
\[
\lim_{n\to\infty} \gamma_n\cdot u \cdot \gamma_n^{-1}
=
\begin{pmatrix}
1 -ac & a^2 \\
- c^2 & 1+ac \\
\end{pmatrix}
\]
In particular, the sequence $\gamma_n\cdot u\cdot \gamma_n^{-1}$
is convergent.
Since $\gamma_n\cdot u\cdot \gamma_n^{-1}\in\Gamma$ for all $n\in\N$,
and
$\Gamma$ is discrete, it follows that
\[
\gamma_n\cdot u \cdot\gamma_n^{-1}
=
\begin{pmatrix}
1 -ac & a^2 \\
- c^2 & 1+ac \\
\end{pmatrix}
\]
for all sufficiently large $n$.
In combination with
\[
\pi\left(\gamma_n\cdot u \cdot\gamma_n^{-1}\right)
=\begin{pmatrix} 1-a_nc_n \\ -c_n^2 \\
\end{pmatrix}
\quad\text{ and }\quad
\pi(\gamma_n)
=\begin{pmatrix} a_n\\ -c_n \\
\end{pmatrix}
\]
it follows that
\[
\pi(\gamma_n)\in\left\{
\begin{pmatrix} a\\ -c \\
\end{pmatrix},
\begin{pmatrix} -a\\ c \\
\end{pmatrix}
\right\}
\]
for all suffiently large $n$.
Because $\lim_{n\to\infty}\pi(\gamma_n)=(a,c)$, this contradicts
our assumption that
$(a_n,c_n)\ne (a,c)$ for infinitely many $n\in\N$.

Thus we have established that $\pi(\Gamma)$ is discrete in $\C^2$.

Now the proof of tameness
follows from Corollary~\ref{overshear-tame}.
\end{proof}

\subsection{Riemann surfaces}

Let $M'$ be a Riemann surface different from $\P_1$ and $\C$ and let
$M=M'\setminus\{p\}$ for some point $p\in M'$. Then $M'$ has the upper
half plane $H^+=\{z\in\C:\Im(z)>0\}$ as universal covering and the
fundamental group $\pi_1(M)$ embedds into $PSL_2(\R)=\Aut(H^+)$
as discrete subgroup. The element of $\pi_1(M)$  corresponding to
a circle around $p$ is embedded into $PSL_2(\R)$ as a unipotent element.
Therefore $\pi_1(M)$ becomes a tame discrete subset in $PSL_2(\C)$.

\subsection{Hyperbolic manifolds}\label{sect-hyper}

A real differentiable manifold $M$ is called (real) hyperbolic
if it admits a complete Riemannian metric with constant negative
sectional curvature.
Such a manifold has the hyperbolic $n$-space $\H_n$ as universal
covering ($n=\dim_\R(M)$). Hence the fundamental group 
$\pi_1(M)$ of $M$ acts by isometries on $\H_n$.  Since the group
of isometries of $\H_3$ is isomorphic to $PSL_2(\C)$, the fundamental
group of a real hyperbolic $3$-fold embedds into $PSL_2(\C)$ as a
discrete subgroup $\Gamma$.
A {\em knot} is a circle $C$  embedded into the $3$-sphere $S^3$. It is know
that the {\em complement} $S^3\setminus C$ is
a hyperbolic $3$-fold with finite volume
unless it is a ``satellite knot'' or a ``torus knot''.
The isotropy group $H$ for the $SL_2(\C)$-action on $\H_3$ is compact.
Therefore $\Gamma$ is a discrete subgroup of finite covolume
in $SL_2(\C)$.
Since the knot complement  $S^3\setminus C$ is obviously
non-compact, it is clear that $SL_2(\C)/\Gamma$ is likewise non-compact.
By the theorem of Kazhdan-Margulis 
(\cite{KM})
it follows that $\Gamma$
contains a unipotent element.
Hence tameness of $\Gamma\simeq\pi_1(M)$
(as a discrete subset of $SL_2(\C$)
follows from Proposition~\ref{SL2-unipotent}.

\section{$\C^*$-factors}

In many results we assumed that we discuss a linear algebraic group
which is {\em character-free}, i.e.,
without non-trivial morphism of algebraic groups with values in $\C^*$.
The background for this assumption is the problem that
$\C^*$ and its powers $(\C^*)^n$ has only few algebraic subgroups,
namely only countably many. In particular there is no continuous family
of projections from $(\C^*)^n$ to quotients by algebraic subgroups.
That probihits the usage of probabilistic methods which we employed
in \cite{JW-TAME1} in order to deduce Proposition~\ref{jw1-maxtorus},
which is crucial
for our instrument of ``threshold sequences''.

To emphasize the problems one encounters in discussing tame discrete
subsets of $(\C^+)^n$ we  point out the proposition below.

\begin{proposition}
Let $G=(\C^*)^n$ with $n\ge 2$. Let $R_k$ be a 
strictly increasing sequence of positive real
numbers. Let $\rho$ be an exhaustion function on $G$.
Then there exists a discrete subset $D$ of $G$ such that
\begin{enumerate}
\item
$\#\{x\in D:\rho(x)\le R_k\}\le k$
\item
For every morphism of complex algebraic groups $F:G\to H$ with
$\dim\ker F>0$ the image $F(D)$ is dense in $F(G)$.
\end{enumerate}
\end{proposition}

\begin{proof}
Without loss of generality we may assume that $F$ is surjective.
The image $F(G)=H$ is connected, commutative and reductive and therefore
isomorphic to some $(\C^*)^k$, $0\le k<n=\dim(G)$.
Any such morphism $F$ is determinant by the induced homomorphism
of fundamental groups, hence there are only countably many such
morphisms $F:(\C^*)^n\to(\C^*)^k$.
It suffices to consider only morphisms to $(\C^*)^{n-1}$, because
every morphism of algebraic groups from $(\C^*)^n$ to $(\C^*)^k$
for $k<n$ fibers through $(\C^*)^{n-1}$.
 
We enumerate all surjective morphism of algebraic
groups $F:(\C^*)^n\to(\C^*)^{n-1}=H$ as $(F_j)_{j\in\N}$.
Now we fix a dense countable subset $S\subset H$.
We fix a bijection $\zeta:\N\to\N\times S$ and choose elements
$g_k\in G$ such that
\begin{enumerate}
\item
$F_j(g_k)=s$ if $\zeta(k)=(j,s)$.
\item
$\rho(g_k)\ge\max\{k,R_k\}$.
\end{enumerate}

Then $D=\{g_k:k\in\N\}$
has the desired properties:
\begin{itemize}
\item
By construction we have $S\subset F_j(D)$ for every $j$.
Therefore $F_j(D)$ is dense in $F_j(G)=H$.
\item
$D$ is discrete,
because $\rho(g_k)\ge k\ \forall k$.
\item
$\#\{x\in D:\rho(x)\le R_k\}\le k$
since $\rho(g_k)\ge R_k\ \forall k$.
\end{itemize}
\end{proof}

\bibliographystyle{line}
\bibliography{tame}

\begin{thebibliography}{10}

\bibitem{AU}
Rafael~B. Andrist and Riccardo Ugolini, `A new notion of tameness', {\em J.
  Math. Anal. App} (1) {\bf 472} (2019), 196--215.

\bibitem{FKZ}
Hubert Flenner, Shulim Kaliman and Mikhail Zaidenberg, `A {G}romov-{W}inkelmann
  type theorem for flexible varieties', {\em J. Eur. Math. Soc. (JEMS)} (11)
  {\bf 18} (2016), 2483--2510.

\bibitem{Forstneric}
Franc Forstneri\v{c}, {\em Stein manifolds and holomorphic mappings}, volume~56
  of {\em Ergebnisse der Mathematik und ihrer Grenzgebiete. 3. Folge. A Series
  of Modern Surveys in Mathematics [Results in Mathematics and Related Areas.
  3rd Series. A Series of Modern Surveys in Mathematics]}, second edition, The
  homotopy principle in complex analysis (Springer, Cham, 2017).

\bibitem{KM}
D.~A. Ka\v{z}dan and G.~A. Margulis, `A proof of {S}elberg's hypothesis', {\em
  Mat. Sb. (N.S.)} {\bf 75 (117)}  (1968), 163--168.

\bibitem{YK}
Yuta Kusakabe, `Elliptic characterization and localization of oka manifolds',
  2018.

\bibitem{Margulis}
G.~A. Margulis, {\em Discrete subgroups of semisimple {L}ie groups}, volume~17
  of {\em Ergebnisse der Mathematik und ihrer Grenzgebiete (3) [Results in
  Mathematics and Related Areas (3)]} (Springer-Verlag, Berlin, 1991).

\bibitem{MM}
Yoz\^{o} Matsushima and Akihiko Morimoto, `Sur certains espaces fibr\'{e}s
  holomorphes sur une vari\'{e}t\'{e} de {S}tein', {\em Bull. Soc. Math.
  France} {\bf 88}  (1960), 137--155.

\bibitem{NW}
Junjiro Noguchi and J\"{o}rg Winkelmann, {\em Nevanlinna theory in several
  complex variables and {D}iophantine approximation}, volume 350 of {\em
  Grundlehren der Mathematischen Wissenschaften [Fundamental Principles of
  Mathematical Sciences]} (Springer, Tokyo, 2014).

\bibitem{RR}
Jean-Pierre Rosay and Walter Rudin, `Holomorphic maps from {${\bf C}^n$} to
  {${\bf C}^n$}', {\em Trans. Amer. Math. Soc.} (1) {\bf 310} (1988), 47--86.

\bibitem{JW-LDS}
J\"org Winkelmann, `Large discrete sets in {S}tein manifolds', {\em Math. Z.}
  {\bf 236}  (2001), 883--901.

\bibitem{JW}
\leavevmode\vrule height 2pt depth -1.6pt width 23pt, `On tameness and growth
  conditions', {\em Doc. Math.} {\bf 13}  (2008), 97--101.

\bibitem{JW-TAME1}
\leavevmode\vrule height 2pt depth -1.6pt width 23pt, `Tame discrete subsets in
  stein manifolds', {\em Journal of the Australian Mathematical Society} {\bf
  ???}  (2018), ???

\end{thebibliography}

\end{document}